\documentclass[12pt,a4paper]{article}
\usepackage{amsfonts}

\usepackage{amssymb}

\usepackage{amsmath}

\newtheorem{theorem}{Theorem}[section]
\newtheorem{definition}[theorem]{Definition}
\newtheorem{lemma}[theorem]{Lemma}

\newtheorem{corollary}[theorem]{Corollary}

\begin{document}
\title{Gr\"{o}bner-Shirshov bases and embeddings of algebras\footnote{Supported by the NNSF of China (No.10771077) and
the NSF of Guangdong Province (No.06025062).} }
\author{
 L.A. Bokut\footnote {Supported by RFBR 01-09-00157, LSS--344.2008.1
and SB RAS Integration grant No. 2009.97 (Russia).} \\
{\small \ School of Mathematical Sciences, South China Normal
University}\\
{\small Guangzhou 510631, P.R. China}\\
{\small Sobolev Institute of Mathematics, Russian Academy of
Sciences}\\
{\small Siberian Branch, Novosibirsk 630090, Russia}\\
{\small  bokut@math.nsc.ru}\\
\\
 Yuqun Chen\footnote {Corresponding author.} \  and Qiuhui Mo\\
{\small \ School of Mathematical Sciences, South China Normal
University}\\
{\small Guangzhou 510631, P.R. China}\\
{\small yqchen@scnu.edu.cn}\\
{\small scnuhuashimomo@126.com}}
 \date{}
\maketitle \noindent\textbf{Abstract:} In this paper, by using
Gr\"{o}bner-Shirshov bases, we show that in the following classes,
each (resp. countably generated) algebra can be embedded into a
simple (resp. two-generated) algebra: associative differential
algebras, associative $\Omega$-algebras, associative
$\lambda$-differential algebras. We show that in the following
classes, each countably generated algebra over a countable field $k$
can be embedded into a simple two-generated algebra: associative
algebras, semigroups, Lie algebras, associative differential
algebras, associative $\Omega$-algebras, associative
$\lambda$-differential algebras. Also we prove that any countably
generated module over a free associative algebra $k\langle X
\rangle$ can be embedded into a cyclic $k\langle X \rangle$-module,
where $|X|>1$. We give another proofs of the well known theorems:
each countably generated group (resp. associative algebra,
semigroup, Lie algebra) can be embedded into a two-generated group
(resp. associative algebra, semigroup, Lie algebra).
 \ \

\noindent \textbf{Key words:} Gr\"{o}bner-Shirshov basis, group,
associative algebra,  Lie algebra, associative differential algebra,
associative $\Omega$-algebra, module.

\noindent {\bf AMS} Mathematics Subject Classification(2000): 20E34,
16Sxx, 17B05, 08Cxx, 16D10, 16S15, 13P10

\section{Introduction}
G. Higman, B.H. Neumann and H. Neuman\cite{HNN49} proved that any
countable group is embeddable into a 2-generated group. It means
that the basic rank of variety of groups is equal to two. In
contrast, for example, the basic ranks of varieties of alternative
and Malcev algebras are equal to infinity (I.P. Shestakov
\cite{She}): there is no such $n$, that a countably generated
alternative (Malcev) algebra can be embeddable into $n$-generated
alternative (Malcev) algebra. Even more, for any $n\geq1$, there
exists an alternative (Malcev) algebra generated by $n+1$ elements
which can not be embedded into an $n$-generated alternative (Malcev)
algebra (V.T. Filippov \cite{Fi1981,Fi1984}). For Jordan algebras,
it is known that the basic rank is bigger than 2, since any
2-generated Jordan algebra is special (A.I. Shirshov \cite{Sh56}),
but there exist (even finitely dimensional) non-special Jordan
algebras (A.A. Albert \cite{Al}). A.I. Malcev \cite{Ma} proved that
any countably generated associative algebra is embeddable into a
2-generated associative algebra. A.I. Shirshov \cite{Sh58} proved
the same result for Lie algebras and T. Evans \cite{Ev} proved the
same result for semigroup.

The first example of finitely generated infinite simple group was
constructed by G. Higman \cite{Hi}. Later P. Hall \cite{Ha} proved
that any group is embeddable into a simple group which is generated
by 3 prescribed subgroups with some cardinality conditions. In
particular, any countably generated group is embeddable into a
simple 3-generated group.

B. Neumann proved that any non-associative algebra is embeddable
into a non-associative division algebra such that any equation
$ax=b$, $xa=b$, $a\neq 0$ has a solution in the latter. Any division
algebra is simple. P.M. Cohn \cite{Co} proved that any associative
ring without zero divisors is embeddable into a simple associative
ring without zero divisors such that any equation $ax - xa=b$,
$a\neq 0$, has a solution in the latter. L.A. Skornyakov \cite{Sk}
proved that any non-associative algebra without zero divisors is
embeddable into a non-associative division algebra without zero
divisors. I.S. Ivanov \cite{Iv65,Iv67} prove the same result for
$\Omega$-algebras (see also A.G. Kurosh \cite{Ku2}). P.M. Cohn
\cite{Co} proved that any Lie algebra is embeddable into a division
Lie algebra. E.G. Shutov \cite{Shu} and L.A. Bokut \cite{Bo63}
proved that any semigroup is embeddable into a simple semigroup, and
L.A. Bokut \cite{Bo76} proved that any associative algebra is
embeddable into a simple associative algebra such that any equation
$xay=b$, $a\neq 0$ is solvable in the latter. L.A. Bokut \cite{
Bo62, Bo62L} proved that any Lie (resp. non-associative,
commutative, anti-commutative) algebra $A$ is embeddable into an
algebraically closed (in particular simple) Lie (resp.
non-associative, commutative, anti-commutative) algebra $B$ such
that any equation $f(x_1,...,x_n)=0$ with coefficient in $B$ has a
solution in $A$ (an equation over $B$ is an element of a free
product of $B$ with a corresponding free algebra $k(X)$). L.A. Bokut
\cite{Bo76, Bo78, BoK} proved that any associative (Lie) algebra is
embeddable into a simple associative (algebraically closed Lie)
algebra which is a sum of 4 prescribed (Lie) subalgebras with some
cardinality conditions. In particular any countable associative
(Lie) algebra is embeddable into a simple finitely generated
associative (Lie) algebra. A.P. Goryushkin \cite{Go} proved that any
countable group is embeddable into a simple 2-generated group.

In this paper, by using Gr\"{o}bner-Shirshov bases and some ideas
from \cite{Bo76, Bo78}, we prove that in the following classes, each
(resp. countably generated) algebra can be embedded into a simple
(resp. two-generated) algebra:  associative differential algebras,
associative $\Omega$-algebras, associative $\lambda$-differential
algebras. We show that in the following classes, each countably
generated algebra over a countable field $k$ can be embedded into a
simple two-generated algebra: associative algebras, semigroups, Lie
algebras, associative differential algebras, associative
$\Omega$-algebras, associative $\lambda$-differential algebras. Also
we prove that any countably generated module over  a free
associative algebra $k\langle X \rangle$ can be embedded into a
cyclic $k\langle X \rangle$-module, where $|X|>1$. We give another
proofs of Higman-Neumann-Neumann's and Shirshov's results mentioned
above.

We systematically use Gr\"{o}bner-Shirshov bases theory for
associative algebras, Lie algebras, associative $\Omega$-algebras,
associative differential algebras, modules, see \cite{Sh62,
 BoCQ, ChCL,ChCZ}.

\section{Preliminaries}

We first cite some concepts and results from the literature
\cite{Bo72, Bo76, Sh62} which are related to Gr\"{o}bner-Shirshov
bases for associative algebras.

Let $X$ be a set and $k$ a field.  Throughout this paper, we denote
by $k\langle X\rangle$ the free associative algebra over $k$
generated by $X$,  by $X^*$ the free monoid generated by $X$ and by
 $N$ the set of natural numbers.

A well ordering $<$ on $X^*$ is called monomial if for $u, v\in
X^*$, we have
$$
u < v \Rightarrow w|_u < w|_v  \ for \  all \
 w\in  X^*,
$$
where $w|_u=w|_{x_i\mapsto u}, \ w|_v=w|_{x_i\mapsto v}$ and
$x_{i}$'s are the same individuality of the letter $x_{i}\in X$ in
$w$.

 A standard example of monomial ordering on $X^*$ is the
deg-lex ordering which first compare two words by degree and then by
comparing them lexicographically, where $X$ is a well ordered set.

Let $X^*$ be a set with a monomial ordering $<$. Then, for any
polynomial $f\in k\langle X\rangle$, $f$ has the leading word
$\overline{f}$. We call $f$  monic if the coefficient of
$\overline{f}$ is 1. By $deg(f)$ we denote the degree of
$\overline{f}$.

Let $f,\ g\in k\langle X\rangle$ be two monic polynomials and $w\in
X^*$. If  $w=\overline{f}b=a\overline{g}$  for some $a,b\in X^*$
such that $deg(\overline{f})+deg(\overline{g})>deg(w)$, then
$(f,g)_w=fb-ag$ is called the intersection composition of $f,g$
relative to $w$. If $w=\overline{f}=a\overline{g}b$  for some $a,
b\in X^*$, then $(f,g)_w=f-agb$ is called the inclusion composition
of $f,g$ relative to $w$. In $(f,g)_w$, $w$ is called the ambiguity
of the composition.

Let $S\subset k\langle X\rangle$ be a monic set. A composition
$(f,g)_w$ is called trivial modulo $(S,w)$, denoted by
$$
(f,g)_w\equiv0 \ \ \ mod(S,w)
$$
if $(f,g)_w=\sum\alpha_ia_is_ib_i,$ where every $\alpha_i\in k, \
s_i\in S,\ a_i,b_i\in X^*$, and $a_i\overline{s_i}b_i<w$.

Recall that $S$ is a Gr\"{o}bner-Shirshov basis in $k\langle
X\rangle$ if any composition of polynomials from $S$ is trivial
modulo $S$.

The following lemma was first proved by Shirshov \cite{Sh62} for
free Lie algebras (with deg-lex ordering) (see also Bokut
\cite{Bo72}). Bokut \cite{Bo76} specialized the approach of Shirshov
to associative algebras (see also Bergman \cite{Be}). For
commutative polynomials, this lemma is known as Buchberger's Theorem
(see \cite{Bu65, Bu70}).

\begin{lemma}\label{l1}
{\em (Composition-Diamond lemma for associative algebras)} \ Let $k$
be a field, $A=k \langle X|S\rangle=k\langle X\rangle/Id(S)$ and $<$
a monomial ordering on $X^*$, where $Id(S)$ is the ideal of $k
\langle X\rangle$ generated by $S$. Then the following statements
are equivalent:

\begin{enumerate}
\item[(1)] $S$ is a Gr\"{o}bner-Shirshov basis in $k\langle
X\rangle$.
\item[(2)] $f\in Id(S)\Rightarrow \bar{f}=a\bar{s}b$
for some $s\in S$ and $a,b\in  X^*$.
\item[(3)] $Irr(S) = \{ u \in X^* |  u \neq a\bar{s}b,s\in S,a ,b \in X^*\}$
is a $k$-basis of the algebra $A=k\langle X | S \rangle=k\langle
X\rangle/Id(S)$.
\end{enumerate}
\end{lemma}

If a subset $S$ of $k\langle X \rangle$ is not a
Gr\"{o}bner-Shirshov basis then one can add all nontrivial
compositions of polynomials of $S$ to $S$. Continue this process
repeatedly, we finally obtain a Gr\"{o}bner-Shirshov basis
$S^{comp}$ that contains $S$. Such a process is called Shirshov
algorithm.

Let $A=sgp\langle X|S\rangle$ be a semigroup presentation. Then $S$
is also a subset of $k\langle X \rangle$ and we can find
Gr\"{o}bner-Shirshov basis
 $S^{comp}$. We also call $S^{comp}$ a
Gr\"{o}bner-Shirshov basis of $A$. $Irr(S^{comp})=\{u\in X^*|u\neq
a\overline{f}b,\ a ,b \in X^*,\ f\in S^{comp}\}$ is a $k$-basis of
$k\langle X|S\rangle$ which is also the set of all normal words of
$A$.

\ \

The following lemma is well known and can be easily proved.

\begin{lemma}\label{l2.2}
Let $k$ be a field, $S\subset k\langle X\rangle$. Then for any $f\in
k\langle X\rangle$, $f$ can be expressed as $f=\sum_{u_i\in
Irr(S),u_i\leq\bar f}\alpha_iu_i +
\sum_{a_j\overline{s_j}b_j\leq\overline{f}}\beta_ja_js_jb_j$, where
$\alpha_i, \beta_j\in k,  a_j, b_j\in X^*, s_j\in S$.
\end{lemma}

The analogous lemma is valid for the free Lie algebra $Lie(X)$ (see,
for example, \cite{BoC1}).

\begin{lemma}\label{l2.3}
Let $k$ be a field, $S\subset Lie(X)$. Then for any $f\in Lie(X)$,
$f$ can be expressed as $f=\sum_{[u_i]\in Irr(S),[u_i]\leq\bar
f}\alpha_i[u_i] +
\sum_{a_j\overline{s_j}b_j\leq\overline{f}}\beta_j[a_js_jb_j]$,
where $\alpha_i, \beta_j\in k,  a_j, b_j\in X^*, s_j\in S$, and
$Irr(S) = \{[u]|[u] \mbox{ is a non-associative Lyndon-Shirshov word
on}\ X, \ u \neq a\bar{s}b,s\in S,a ,b \in X^*\}$.
\end{lemma}

 \ \

\section{Associative algebras, Groups and Semigroups}
In this section we give another proofs for the following theorems
mentioned in the introduction: every countably generated group
(resp. associative algebra, semigroup) can be embedded into a
two-generated group (resp. associative algebra, semigroup). Even
more, we prove the following theorems: (i) Every countably generated
associative algebra over a countable field $k$ can be embedded into
a simple two-generated associative algebra. (ii) Every countably
generated semigroup can be embedded into a (0-)simple two-generated
semigroup.

\ \

In this section, all the algebras we mention contain units.

\ \

In 1949, G. Higman, B.H. Neumann, and H. Neumann \cite{HNN49} prove
that every countable group can be embedded into a two-generated
group. Now we give another proof for this theorem.
\begin{theorem}
(G. Higman, B.H. Neumann and H. Neumann) \  Every countable group
can be embedded into a two-generated group.
\end{theorem}
\textbf{Proof}\ We may assume that the group
$G=\{g_{0}=1,g_{1},g_{2},g_{3},\dots \}$. Let
$$
H=gp\langle G\setminus\{g_0\},a,b,t|g_{j}g_{k}=[g_{j},g_{k}],\
at=tb, b^{-i}ab^{i}t=tg_{i} a^{-i}ba^{i},\ i,j,k\in N \rangle.
$$
G. Higman, B.H. Neumann and H. Neumann \cite{HNN49} (see also
\cite{LySc}) proved that $G$ can be embedded into $H$. Now, we use
the Composition-Diamond lemma, i.e., Lemma \ref{l1} to reprove this
theorem.

Clearly, $H$ can also be expressed as
$$
H=gp\langle G\setminus\{g_0\},a,b,t| S \rangle,
$$
where $S$ consists of the following polynomials ($\varepsilon=\pm
1,i,j,k\in N$):
\begin{eqnarray*}
1.&&g_{j}g_{k}=[g_{j},g_{k}]\\
2.&&a^\varepsilon t=tb^{\varepsilon}\\
3.&&b^{\varepsilon} t^{-1}=t^{-1}a^{\varepsilon}\\
4.&&ab^{i}t=b^{i}tg_{i}a^{-i}ba^{i}\\
5.&&a^{-1}b^{i}t=b^{i}t(g_{i}a^{-i}ba^{i})^{-1}\\
6.&&ba^{i}t^{-1}=a^{i}g_{i}^{-1}t^{-1}b^{-i}ab^{i}\\
7.&&b^{-1}a^{i}g_{i}^{-1}t^{-1}=a^{i}t^{-1}b^{-i}a^{-1}b^{i}\\
8.&&a^{\varepsilon} a^{-\varepsilon}=b^{\varepsilon}
b^{-\varepsilon}=t^{\varepsilon} t^{-\varepsilon}=1
\end{eqnarray*}
We order $\{g_{i},a^{\pm 1}, b^{\pm 1}\}^{*} $ by deg-lex ordering
with $ g_{i}<a<a^{-1}<b<b^{-1}. $ Denote by $X=\{g_{i},a^{\pm 1},
b^{\pm 1}, t^{\pm 1}\}$. For any $u \in X^{*}$, $u$ can be uniquely
expressed without brackets as
$$
u=u_{0}t^{\varepsilon_{1}}u_{1}t^{\varepsilon_{2}}u_{2}\cdots
t^{\varepsilon_{n}}u_{n},
$$
where $u_{i}\in \{g_{i},a^{\pm 1}, b^{\pm 1}\}^{*} , n\geq 0,
\varepsilon_{i}=\pm 1$. Denote by
$$
wt(u)=(n,u_{0},t^{\varepsilon_{1}},u_{1},t^{\varepsilon_{2}},u_{2},\dots,
t^{\varepsilon_{n}},u_{n}).
$$
Then, we order $X^*$ as follows: for any $u,v\in X^*$
$$
u>v\Leftrightarrow wt(u)>wt(v)\ \ \ \ \ \ \ \  lexicographically,
$$
where $t>t^{-1}$. With this ordering, we can check that $S$ is a
Gr\"{o}bner-Shirshov basis in the free associative algebra $k\langle
X\rangle$. By Lemma \ref{l1}, the group $G$ can be embedded into $H$
which is generated by $\{a,b\}$. \hfill $\blacksquare$

\ \

A.I. Malcev \cite{Ma} proved that any countably generated
associative algebra  is embeddable into a two-generated associative
algebra, and T. Evans \cite{Ev} proved that every countably
generated semigroup  can be embedded into a two-generated semigroup.
Now, by applying Lemma \ref{l1}, we give another proofs of this two
embedding theorems.

\begin{theorem}\label{t4.2}
(A.I. Malcev) \ Every countably generated associative algebra can be
embedded into a two-generated associative algebra.
\end{theorem}
\textbf{Proof} \ We may assume that $ A= k \langle X  |S\rangle $ is
an associative algebra  generated by $X$ with relations $S$, where
$X=\{x_{i},i=1,2,\dots\}$. By Shirshov algorithm, we can assume that
$S$ is a Gr\"{o}bner-Shirshov basis in the free associative algebra
$k\langle X\rangle$ with deg-lex ordering on $X^*$.  Let
$$
H=k\langle X, a,b | S, aab^{i}ab=x_{i}, \ i=1,2,\dots\rangle.
$$
We can  check that
$$\{S,
aab^{i}ab=x_{i}, \ i=1,2,\dots\}
$$
is a Gr\"{o}bner-Shirshov basis in  $ k\langle X,a,b \rangle$ with
deg-lex ordering on $(X\cup\{a,b\})^*$ where $a>b>x,\ x\in X$ since
there are no new compositions. By Lemma \ref{l1}, $A$ can be
embedded into $H$ which is generated by $\{a,b\}$. \hfill
$\blacksquare$

\ \

By the proof of Theorem \ref{t4.2}, we have immediately the
following corollary.
\begin{corollary}
(T. Evans) \  Every countably generated semigroup can be embedded
into a two-generated semigroup.
\end{corollary}

\begin{theorem}\label{t0.0}
Every countably generated associative algebra over a countable field
$k$ can be embedded into a simple two-generated associative algebra.
\end{theorem}
\textbf{Proof} \ Let $A$ be a  countably generated associative
algebra over a countable field $k$. We may assume that $A$ has a
countable $k$-basis $\{1\}\cup X_0$, where $X_0=\{x_{i}|
i=1,2,\ldots\}$ and 1 is the unit of $A$. Then $A$ can be expressed
as $A= k\langle X_0 |x_{i}x_{j}=\{x_{i},x_{j}\}, i,j\in N \rangle$,
where $\{x_{i},x_{j}\}$ is a linear combination of $x_t, \ x_t\in
X_0$.

Let $A_0=k\langle X_0\rangle$, $A_0^+=A_0\backslash\{0\}$ and fix
the bijection
$$
(A_0^+,A_0^+)\longleftrightarrow\{(x_m^{(1)},y_m^{(1)}), m\in N\}.
$$

Let $X_1=X_0\cup\{x_m^{(1)},y_m^{(1)},a,b|m\in N\}$, $A_1=k\langle
X_1\rangle$, $A_1^+=A_1\backslash\{0\}$ and fix the bijection
$$
(A_1^+,A_1^+)\longleftrightarrow\{(x_m^{(2)},y_m^{(2)}), m\in N\}.
$$
$$
\vdots
$$
Let $X_{n+1}=X_n\cup\{x_m^{(n+1)},y_m^{(n+1)}|m\in N\}$, $n\geq1$,
$A_{n+1}=k\langle X_{n+1}\rangle$,
$A_{n+1}^+=A_{n+1}\backslash\{0\}$ and fix the bijection
$$
(A_{n+1}^+,A_{n+1}^+)\longleftrightarrow\{(x_m^{(n+2)},y_m^{(n+2)}),
m\in N\}.
$$
$$
\vdots
$$

Consider the chain of the free associative algebras
$$
A_0\subset A_1\subset A_2\subset\ldots\subset A_n\subset\ldots.
$$

Let $X=\bigcup_ {n=0}^ {\infty}X_n$. Then $k\langle
X\rangle=\bigcup_ {n=0}^ {\infty}A_n$.

Now, define the desired algebra $\mathcal{A}$. Take the set $X$ as
the set of the generators for this algebra and take the following
relations as one part of the relations for this algebra
\begin{equation}\label{1}
x_{i}x_{j}=\{x_{i},x_{j}\}, \ i,j\in N
\end{equation}
\begin{equation}\label{2}
aa(ab)^nb^{2m+1}ab=x_m^{(n)},\  m,n\in N
\end{equation}
\begin{equation}\label{3}
aa(ab)^nb^{2m}ab=y_m^{(n)}, \ m,n\in N
\end{equation}
\begin{equation}\label{4}
aabbab=x_1
\end{equation}

Before we introduce the another part of the relations on
$\mathcal{A}$, let us define canonical words of the algebra $A_n$,
$n\geq0$. A word in $X_0$ without subwords that are the leading
terms of (\ref{1}) is called a canonical word of $A_0$. A word in
$X_1$ without subwords that are the leading terms of (\ref{1}),
(\ref{2}), (\ref{3}), (\ref{4}) and without subwords of the form
$$
(x_m^{(1)})^{deg(g^{(0)})+1}\overline{f^{(0)}}y_m^{(1)},
$$
where
$(x_m^{(1)},y_m^{(1)})\longleftrightarrow(f^{(0)},g^{(0)})\in(A_0^+,A_0^+)$
such that $f^{(0)},g^{(0)}$ are non-zero linear combination of
canonical words of $A_0$, is called a canonical word of $A_1$.
Suppose that we have defined canonical word of $A_{k}$, $k<n$. A
word in $X_n$ without subwords that are the leading terms of
(\ref{1}), (\ref{2}), (\ref{3}), (\ref{4}) and without subwords of
the form
$$
(x_m^{(k+1)})^{deg(g^{(k)})+1}\overline{f^{(k)}}y_m^{(k+1)},
$$
where
$(x_m^{(k+1)},y_m^{(k+1)})\longleftrightarrow(f^{(k)},g^{(k)})\in(A_{k}^+,A_{k}^+)$
such that  $f^{(k)},g^{(k)}$ are non-zero linear combination of
canonical words of $A_{k}$, is called a canonical word of $A_{n}$.

Then the another part of the relations on $\mathcal{A}$ are the
following:
\begin{equation}\label{5}
(x_m^{(n)})^{deg(g^{(n-1)})+1}f^{(n-1)}y_m^{(n)}-g^{(n-1)}=0,\
m,n\in N
\end{equation}
where
$(x_m^{(n)},y_m^{(n)})\longleftrightarrow(f^{(n-1)},g^{(n-1)})\in(A_{n-1}^+,A_{n-1}^+)$
such that $f^{(n-1)},g^{(n-1)}$ are non-zero linear combination of
canonical words of $A_{n-1}$.

By Lemma \ref{l2.2}, we have that in $\mathcal{A}$ every element can
be expressed as linear combination of canonical words.

Denote by $S$ the set constituted by the relations
(\ref{1})-(\ref{5}). We can see that $S$ is a Gr\"{o}bner-Shirshov
basis in $k \langle X\rangle$ with deg-lex ordering on $X^*$ since
in $S$ there are no compositions except for the ambiguity
$x_{i}x_{j}x_{k}$ which is a trivial case.  By Lemma \ref{l1}, $A$
can be embedded into $\mathcal{A}$. By (\ref{5}), $\mathcal{A}$ is a
simple algebra. By (\ref{2})-(\ref{5}), $\mathcal{A}$ is generated
by $\{a,b\}$. \hfill $\blacksquare$

\ \

A semigroup $S$ without zero is called simple if it has no proper
ideals. A semigroup $S$ with zero is called 0-simple if $\{0\}$ and
$S$ are its only ideals, and $S^2\neq\{0\}$.

\begin{lemma}(\cite{Ho})
A semigroup $S$ with 0 is 0-simple if and only if $SaS=S$ for every
$a\neq0$ in $S$. A semigroup $S$ without 0 is simple if and only if
$SaS=S$ for every $a$ in $S$.
\end{lemma}

 The following theorem follows from the proof of Theorem \ref{t0.0}.

\begin{theorem}\label{t3.0}
Every countably generated semigroup can be embedded into a simple
two-generated semigroup.
\end{theorem}

\ \

\noindent\textbf{Remark:} Let ${S}$ be a simple semigroup. Then the
semigroup ${S}^0$ with 0 attached is a 0-simple semigroup.
Therefore, by Theorem \ref{t3.0}, each countably generated semigroup
can be embedded into a 0-simple two-generated semigroup.

\ \

\section{Lie algebras}

In this section, we give another proof of the Shirshov's theorem:
every countably generated Lie algebra can be embedded into a
two-generated Lie algebra. And even more, we show that every
countably generated Lie algebra over a countable field $k$ can be
embedded into a simple two-generated Lie algebra.

We start with the Lyndon-Shirshov associative words.

Let $X=\{x_i|i\in I\}$ be a well-ordered set with $x_i>x_p$ if $i>p$
for any $i,p\in I$.  We order $X^*$ by the lexicographical ordering.

\begin{definition}(\cite{Ly,Sh58}, see \cite{BoC1,Uf})
Let $u\in X^\ast$ and $u\neq 1$. Then $u$ is called an $ALSW$
(associative Lyndon-Shirshov word) if
$$
(\forall v,w\in X^\ast, v,w\neq 1) \ u=vw\Rightarrow vw>wv.
$$
\end{definition}

\begin{definition}(\cite{Chen,Sh58}, see \cite{BoC1,Uf})
A non-associative word  $(u)$ in $X$ is called a $NLSW$
(non-associative Lyndon-Shirshov word)  if
\begin{enumerate}
\item[(i)] $u$ is an $ALSW$,
\item[(ii)] if $(u)=((v)(w))$, then both $(v)$ and $(w)$ are
$NLSW$'s,
\item[(iii)] in (ii) if $(v)=((v_1)(v_2))$, then $v_2\leq w$ in
$X^\ast$.
\end{enumerate}
\end{definition}

\begin{lemma}(\cite{Chen,Sh58}, see \cite{BoC1,Uf})
Let $u$ be an $ALSW$. Then there exists a unique bracketing way such
that $(u)$ is a $NLSW$.
\end{lemma}

Let $X^{\ast\ast}$ be the set of all non-associative words $(u)$ in
$X$. If $(u)$ is a $NLSW$, then we denote it by $[u]$.

\begin{lemma}(\cite{Chen,Sh58}, see \cite{BoC1,Uf})
$NLSW$'s forms a linear basis of $Lie(X)$,  the free Lie algebra
generated by $X$.
\end{lemma}

Composition-Diamond lemma for free Lie algebras (with deg-lex
ordering) is given in \cite{Sh62} (see also \cite{BoC1}). By
applying this lemma, we give the following theorem.

\begin{theorem}
(A.I. Shirshov) \ Every countably generated Lie algebra can be
embedded into a two-generated Lie algebra.
\end{theorem}
\textbf{Proof}\ We may assume that
$$
L=Lie( X|S)
$$
is a Lie algebra generated by $X$ with relations $S$, where
$X=\{x_{i},i=1,2,\dots\}$. By Shirshov algorithm, we can assume that
$S$ is a Gr\"{o}bner-Shirshov basis  in the free Lie algebra
$Lie(X)$ on deg-lex ordering. Let
$$
H=Lie( X, a,b | S, [aab^{i}ab]=x_{i}, \ i=1,2,\dots).
$$
We can  check that
$$
\{ S, [aab^{i}ab]=x_{i}, \ i=1,2,\dots\}
$$
is a Gr\"{o}bner-Shirshov basis of $Lie( X,a,b )$ on deg-lex
ordering with $a>b>x_i$ since there are no new compositions. By the
Composition-Diamond lemma for Lie algebras, $L$ can be embedded into
$H$ which is generated by $\{a,b\}$. \hfill $\blacksquare$

\begin{theorem}
Every countably generated Lie algebra over a countable field $k$ can
be embedded into a simple two-generated Lie algebra.
\end{theorem}
\textbf{Proof} \ Let $L$ be a countably generated Lie algebra over a
countable field $k$. We may assume that $L$ has a countable
$k$-basis  $X_0=\{x_{i}| i=1,2,\ldots\}$. Then $L$ can be expressed
as $L= Lie(X_0 |[x_{i}x_{j}]=\{x_{i},x_{j}\}, i,j\in N)$.

Let $L_0=Lie(X_0)$, $L_0^+=L_0\backslash\{0\}$ and fix the bijection
$$
(L_0^+,L_0^+)\longleftrightarrow\{(x_m^{(1)},y_m^{(1)}), m\in N\}.
$$

Let $X_1=X_0\cup\{x_m^{(1)},y_m^{(1)},a,b|m\in N\}$, $L_1=Lie(
X_1)$, $L_1^+=L_1\backslash\{0\}$ and fix the bijection
$$
(L_1^+,L_1^+)\longleftrightarrow\{(x_m^{(2)},y_m^{(2)}), m\in N\}.
$$
$$
\vdots
$$
Let $X_{n+1}=X_n\cup\{x_m^{(n+1)},y_m^{(n+1)}|m\in N\}$, $n\geq1$,
$L_{n+1}=Lie(X_{n+1})$, $L_{n+1}^+=L_{n+1}\backslash\{0\}$ and fix
the bijection
$$
(L_{n+1}^+,L_{n+1}^+)\longleftrightarrow\{(x_m^{(n+2)},y_m^{(n+2)}),
m\in N\}.
$$
$$
\vdots
$$

Consider the chain of the free Lie algebras
$$
L_0\subset L_1\subset L_2\subset\ldots\subset L_n\subset\ldots.
$$

Let $X=\bigcup_ {n=0}^ {\infty}X_n$. Then $Lie(X)=\bigcup_ {n=0}^
{\infty}L_n$.

Now, define the desired Lie algebra $\mathcal{L}$. Take the set $X$
as the set of the generators for this algebra and take the following
relations as one part of the relations for this algebra
\begin{equation}\label{11}
[x_{i}x_{j}]=\{x_{i},x_{j}\}, \ i,j\in N
\end{equation}
\begin{equation}\label{12}
[aa(ab)^nb^{2m+1}ab]=x_m^{(n)},\  m,n\in N
\end{equation}
\begin{equation}\label{13}
[aa(ab)^nb^{2m}ab]=y_m^{(n)}, \ m,n\in N
\end{equation}
\begin{equation}\label{14}
[aabbab]=x_1
\end{equation}

Before we introduce the another part of the relations on
$\mathcal{L}$, let us define canonical words of the Lie algebra
$L_n$, $n\geq0$. A NLSW $[u]$ in $X_0$ where $u$ without subwords
that are the leading terms of (\ref{11}) is called a canonical word
of $L_0$. A NLSW $[u]$ in $X_1$ where $u$ without subwords that are
the leading terms of (\ref{11}), (\ref{12}), (\ref{13}), (\ref{14})
and without subwords of the form
$$
x_m^{(1)}\overline{f^{(0)}}x_m^{(1)}(y_m^{(1)})^{deg(g^{(0)})+1},
$$
where
$(x_m^{(1)},y_m^{(1)})\longleftrightarrow(f^{(0)},g^{(0)})\in(L_0^+,L_0^+)$
such that $f^{(0)},g^{(0)}$ are non-zero linear combination of
canonical words of $L_0$, is called a canonical word of $L_1$.
Suppose that we have defined canonical word of $L_{k}$, $k<n$. A
NLSW $[u]$ in $X_n$ where $u$ without subwords that are the leading
terms of (\ref{11}), (\ref{12}), (\ref{13}), (\ref{14}) and without
subwords of the form
$$
x_m^{(k+1)}\overline{f^{(k)}}x_m^{(k+1)}(y_m^{(k+1)})^{deg(g^{(k)})+1},
$$
where
$(x_m^{(k+1)},y_m^{(k+1)})\longleftrightarrow(f^{(k)},g^{(k)})\in(L_k^+,L_k^+)$
such that $f^{(k)},g^{(k)}$ are non-zero linear combination of
canonical words of $L_k$, is called a canonical word of $L_n$.

Then the another part of the relations on $\mathcal{L}$ are the
following:
\begin{equation}\label{15}
(x_m^{(n)}f^{(n-1)})[x_m^{(n)}(y_m^{(n)})^{deg(g^{(n-1)})+1}]-g^{(n-1)}=0,\
\ m, n\in N
\end{equation}
where
$(x_m^{(n)},y_m^{(n)})\longleftrightarrow(f^{(n-1)},g^{(n-1)})\in(L_{n-1}^+,L_{n-1}^+)$
such that  $f^{(n-1)},g^{(n-1)}$ are non-zero linear combination of
canonical words of $L_{n-1}$.

By Lemma \ref{l2.3}, we have in $\mathcal{L}$ every element can be
expressed as linear combination of canonical words.

Denote by $S$ the set constituted by the relations
(\ref{11})-(\ref{15}).  Define
$\ldots>x_q^{(2)}>x_m^{(1)}>a>b>x_i>y_n^{(1)}>y_p^{(2)}>\ldots$. We
can see that in $S$ there are no compositions unless for  the
ambiguity $x_{i}x_{j}x_{k}$. But this case is trivial. Hence $S$ is
a Gr\"{o}bner-Shirshov basis in $Lie(X)$ on deg-lex ordering which
implies that $L$ can be embedded into $\mathcal{L}$. By
(\ref{12})-(\ref{15}), $\mathcal{L}$ is a simple Lie algebra
generated by $\{a,b\}$.  \hfill $\blacksquare$

\section{Associative differential algebras}

Composition-Diamond lemma for associative differential algebras with
unit is established in a recent paper \cite{ChCL}. By applying this
lemma in this section, we show that: (i). Every countably generated
associative differential algebra can be embedded into a
two-generated associative differential algebra. (ii). Any
associative differential algebra can be embedded into a simple
associative differential algebra. (iii). Every countably generated
associative differential algebra with countable set $\mathcal{D}$ of
differential operations over a countable field $k$ can be embedded
into a simple two-generated associative differential algebra.

Let $\mathcal{A}$ be an associative algebra over a field $k$ with
unit. Let $\mathcal{D}$ be a set of  multiple linear operations on
$\mathcal{A}$. Then $\mathcal{A}$ is called an associative
differential algebra with differential operations $\mathcal{D}$ or
$\mathcal{D}$-algebra, for short, if for any $D\in \mathcal{D}, \ a,
b \in \mathcal{A}$,
$$
D(ab)=D(a)b+aD(b).
$$

Let $\mathcal{D}=\{D_j|j\in J\}$. For any $m=0,1,\cdots$ and
$\bar{j}=(j_1,\cdots,j_m)\in J^m$, denote by
$D^{\bar{j}}=D_{j_1}D_{j_2}\cdots D_{j_{m}}$ and
$D^{\omega}(X)=\{D^{\bar{j}}(x)|x\in X, \ \bar{j}\in J^m, \ m\geq
0\}$, where $D^{0}(x)=x$.  Let $T=(D^{\omega}(X))^*$ be the free
monoid generated by $D^{\omega}(X)$. For any
$u=D^{\overline{i_1}}(x_1)D^{\overline{i_2}}(x_2)\cdots
D^{\overline{i_n}}(x_n)$ $\in T$, the length of $u$, denoted by
$|u|$, is defined to be $n$. In particular, $|1|=0$.

Let $k\langle X;\mathcal{D} \rangle=kT$ be the $k$-algebra spanned
by $T$. For any $D_j\in \mathcal{D}$,
 we define the linear map $ D_j: \ k\langle X;\mathcal{D} \rangle \rightarrow
k\langle X;\mathcal{D} \rangle$ by induction on $|u|$ for $u\in T$:
\begin{enumerate}
\item[1)] $D_j(1)=0$.
\item[2)] Suppose that $u=D^{\bar{i}}(x)=D_{i_1}D_{i_2}\cdots
D_{i_{m}}(x)$. Then $D_j(u)=D_jD_{i_1}D_{i_2}\cdots D_{i_{m}}(x)$.
\item[3)] Suppose that $u=D^{\bar{i}}(x)\cdot v, \ v\in T$. Then
$D_j(u)=(D_jD^{\bar{i}}(x))\cdot v+D^{\bar{i}}(x)\cdot D_j(v)$.
\end{enumerate}
Then, $k\langle X;\mathcal{D} \rangle$ is  a free associative
differential algebra generated by $X$ with differential operators
$\mathcal{D}$ (see \cite{ChCL}).

Let $\mathcal{D}=\{D_j|j\in J\}$, $X$ and $J$ well ordered sets,
 $D^{\bar{i}}(x)=D_{i_1}D_{i_2}\cdots D_{i_{m}}(x)\in
 D^{\omega}(X)$. Denote by
$$
wt(D^{\bar{i}}(x))=(x; m, i_1, i_2,\cdots, i_{m}).
$$
Then, we order $D^{\omega}(X)$ as follows:
$$
D^{\bar{i}}(x)< D^{\bar{j}}(y)\Longleftrightarrow
wt(D^{\bar{i}}(x))< wt(D^{\bar{j}}(y)) \ \mbox{ lexicographically}.
$$
It is easy to check this ordering is a well ordering on
$D^{\omega}(X)$.

Now, we order $T=(D^{\omega}(X))^*$ by deg-lex ordering. We will use
this ordering in this section.

For convenience, for any $u\in T$, we denote
$\overline{D^{\bar{j}}(u)}$ by ${d^{\bar{j}}(u)}$.

\begin{theorem}
Every countably generated associative differential algebra  can be
embedded into a two-generated associative differential algebra.
\end{theorem}

\textbf{Proof}\  Suppose that $\mathcal{A}=k\langle X;\mathcal{D}|S
\rangle$ is an  associative differential algebra  generated by $X$
with relations $S$, where $X=\{x_{i},i=1,2,\dots\}$. By Shirshov
algorithm, we can assume that with the deg-lex ordering on
$(D^{\omega}(X))^*$ defined as above, $S$ is a Gr\"{o}bner-Shirshov
basis of the free associative differential algebra $k\langle
X;\mathcal{D}\rangle$ in the sense of the paper \cite{ChCL}. Let
$\mathcal{B}=k\langle X,a,b; \mathcal{D}| S, aab^{i}ab=x_{i}
\rangle$. We have that with the deg-lex ordering on
$(D^{\omega}(X,a,b))^*$, $\{S, aab^{i}ab=x_{i}, \ i=1,2,\dots\}$ is
a Gr\"{o}bner-Shirshov basis in the free associative differential
algebra $k\langle X,a,b;\mathcal{D}\rangle$ since there are no new
compositions. By the Composition-Diamond lemma in \cite{ChCL},
$\mathcal{A}$ can be embedded into $\mathcal{B}$ which is generated
by $\{a,b\}$. \hfill $\blacksquare$

\begin{theorem}
Every  associative differential algebra  can be embedded into a
simple associative differential algebra.
\end{theorem}

\textbf{Proof}\ Let ${A}$ be an associative differential algebra
over a field $k$ with $k$-basis $\{1\}\cup X$, where $X=\{x_i\mid
i\in I\}$ and $I$ is a well ordered set.

It is clear that $S_0=\{x_{i}x_{j}=\{x_{i},x_{j}\}, \
D(x_i)=\{D(x_i)\},  i,j\in I,\ D\in \cal D\}$ where $\{D(x_i)\}$ is
a linear combination of $x_j, j\in I$, is a Gr\"{o}bner-Shirshov
basis in the free associative differential algebra $k\langle
X;\mathcal{D}\rangle$ with the deg-lex ordering on
$(D^{\omega}(X))^*$, and ${A}$ can be expressed as
$$
{A}=k\langle X;\mathcal{D}|x_ix_j=\{x_i,x_j\}, D(x_i)=\{D(x_i)\},
 i,j\in I,\ D\in \cal D\rangle.
$$

Let us totally order the set of monic elements of ${A}$. Denote by
$T$ the set of indices for the resulting totally ordered set.
Consider the totally ordered set
$T^2=\{(\theta,\sigma)|\theta,\sigma\in T\}$ and assign
$(\theta,\sigma)<(\theta',\sigma')$ if either $\theta<\theta'$ or
$\theta=\theta'$ and $\sigma<\sigma'$. Then $T^2$ is also totally
ordered set.

For each ordered pair of elements $f_\theta, f_\sigma\in{A}, \ \
\theta, \sigma\in T$, introduce the letters
$x_{\theta\sigma},y_{\theta\sigma}$.

Let ${A}_{1}$ be the associative differential algebra given by the
generators
$$
X_1=\{x_i,y_{\theta\sigma},x_{\varrho\tau}| i\in I,\ \theta, \sigma,
\varrho, \tau\in T\}
$$
and the defining relations
$$
S=\{x_ix_j=\{x_i,x_j\}, D(x_i)=\{D(x_i)\}, x_{\theta\sigma}f_\theta
y_{\theta\sigma}=f_\sigma\mid i, j\in I, \ D\in {\cal D},\
({\theta,\sigma})\in T^2\}.
$$
We can have that with the deg-lex ordering on $(D^{\omega}(X_1))^*$,
$S$ is a Gr\"{o}bner-Shirshov basis in the free associative
differential algebra $k\langle X_1;\mathcal{D}\rangle$ in the sense
of the paper \cite{ChCL} since there are no new compositions. Thus,
by the Composition-Diamond lemma in \cite{ChCL}, ${A}$ can be
embedded into ${A}_1$. The relations $x_{\theta\sigma}f_\theta
y_{\theta\sigma}=f_\sigma$ of ${A}_1$ provide that in ${A}_1$ every
monic element $f_\theta$ of the subalgebra ${A}$ generates an ideal
containing algebra $A$.

Mimicking the construction of the associative differential algebra
${A}_1$ from the ${A}$, produce the associative differential algebra
${A}_2$ from ${A}_1$ and so on. As a result, we acquire an ascending
chain of associative differential algebras
$$
{A}={A}_0\subset{A}_1\subset{A}_2\subset\cdots
$$
such that every  monic element $f\in A_k$ generates an ideal in
${A}_{k+1}$ containing $A_k$. Therefore, in the associative
differential algebra
$$
{\cal A}=\bigcup_ {k=0}^ {\infty}{A}_k,
$$
every  nonzero element generates the same ideal. Thus, ${\cal A}$ is
a simple associative differential algebra. \hfill $\blacksquare$

\begin{theorem}
Every countably generated associative differential algebra with
countable set $\mathcal{D}$ of differential operations over a
countable field $k$ can be embedded into a simple two-generated
associative differential algebra.
\end{theorem}
\textbf{Proof} \ Let $A$ be a countably generated associative
differential algebra with countable set $\mathcal{D}$ of
differential operations over a countable field $k$. We may assume
that $A$ has a countable $k$-basis $\{1\}\cup X_0$, where
$X_0=\{x_{i}| i=1,2,\ldots\}$. Then $A$ can be expressed as
$$
A=k\langle X_0;\mathcal{D}|x_ix_j=\{x_i,x_j\}, D(x_i)=\{D(x_i)\},
i,j\in N,\ D\in \cal D\rangle.
$$

Let $A_0=k\langle X_0;\mathcal{D}\rangle$,
$A_0^+=A_0\backslash\{0\}$ and fix the bijection
$$
(A_0^+,A_0^+)\longleftrightarrow\{(x_m^{(1)},y_m^{(1)}), m\in N\}.
$$

Let $X_1=X_0\cup\{x_m^{(1)},y_m^{(1)},a,b|m\in N\}$, $A_1=k\langle
X_1;\mathcal{D}\rangle$, $A_1^+=A_1\backslash\{0\}$ and fix the
bijection
$$
(A_1^+,A_1^+)\longleftrightarrow\{(x_m^{(2)},y_m^{(2)}), m\in N\}.
$$
$$
\vdots
$$

Let $X_{n+1}=X_n\cup\{x_m^{(n+1)},y_m^{(n+1)}|m\in N\}$, $n\geq1$,
$A_{n+1}=k\langle X_{n+1};\mathcal{D}\rangle$,
$A_{n+1}^+=A_{n+1}\backslash\{0\}$ and fix the bijection
$$
(A_{n+1}^+,A_{n+1}^+)\longleftrightarrow\{(x_m^{(n+2)},y_m^{(n+2)}),
m\in N\}.
$$
$$
\vdots
$$

Consider the chain of the free associative differential algebras
$$
A_0\subset A_1\subset A_2\subset\ldots\subset A_n\subset\ldots.
$$

Let $X=\bigcup_{n=0}^{\infty}X_n$. Then $k\langle
X;\mathcal{D}\rangle=\bigcup_{n=0}^{\infty}A_n$.

Now, define the desired associative differential algebra
$\mathcal{A}$. Take the set $X$ as the set of the generators for
this algebra and take the following relations as one part of the
relations for this algebra
\begin{equation}\label{16}
x_{i}x_{j}=\{x_{i},x_{j}\}, \ D(x_i)=\{D(x_i)\}, \ i,j\in N,\ D\in
\cal D
\end{equation}
\begin{equation}\label{17}
aa(ab)^nb^{2m+1}ab=x_m^{(n)},\  m,n\in N
\end{equation}
\begin{equation}\label{18}
aa(ab)^nb^{2m}ab=y_m^{(n)}, \ m,n\in N
\end{equation}
\begin{equation}\label{19}
aabbab=x_1
\end{equation}

Before we introduce the another part of the relations on
$\mathcal{A}$, let us define canonical words of the algebras $A_n$,
$n\geq0$. An element in $(D^{\omega}(X_0))^*$ without subwords of
the form ${d^{\bar{i}}(u)}$ where $u$ is the leading terms of
(\ref{16}), is called a canonical word of $A_0$. An element in
$(D^{\omega}(X_1))^*$ without subwords of the form
${d^{\bar{i}}(u)}$ where $u$ is the leading terms of (\ref{16}),
(\ref{17}), (\ref{18}), (\ref{19}) and
$$
(x_m^{(1)})^{|\overline{g^{(0)}}|+1}\overline{f^{(0)}}y_m^{(1)},
$$
where
$(x_m^{(1)},y_m^{(1)})\longleftrightarrow(f^{(0)},g^{(0)})\in(A_0^+,A_0^+)$
such that $f^{(0)},g^{(0)}$ are non-zero linear combination of
canonical words of $A_0$, is called a canonical word of $A_1$.
Suppose that we have defined canonical word of $A_{k}$, $k<n$. An
element in $(D^{\omega}(X_n))^*$ without subwords of the form
${d^{\bar{i}}(u)}$ where $u$ is the leading terms of (\ref{16}),
(\ref{17}), (\ref{18}), (\ref{19}) and
$$
(x_m^{(k+1)})^{|\overline{g^{(k)}}|+1}\overline{f^{(k)}}y_m^{(k+1)},
$$
where
$(x_m^{(k+1)},y_m^{(k+1)})\longleftrightarrow(f^{(k)},g^{(k)})\in(A_k^+,A_k^+)$
such that $f^{(k)},g^{(k)}$ are non-zero linear combination of
canonical words of $A_k$, is called a canonical word of $A_n$.

Then the another part of the relations on $\mathcal{A}$ are the
following:
\begin{equation}\label{20}
(x_m^{(n)})^{|\overline{g^{(n-1)}}|+1}f^{(n-1)}y_m^{(n)}-g^{(n-1)}=0,
\ \ m, n\in N
\end{equation}
where
$(x_m^{(n)},y_m^{(n)})\longleftrightarrow(f^{(n-1)},g^{(n-1)})\in(A_{n-1}^+,A_{n-1}^+)$
such that $f^{(n-1)},g^{(n-1)}$ are non-zero linear combination of
canonical words of $A_{n-1}$.

We can get that in $\mathcal{A}$ every element can be expressed as
linear combination of canonical words.

Denote by $S$ the set constituted by the relations
(\ref{16})-(\ref{20}). We can have that with the deg-lex ordering on
$(D^{\omega}(X))^*$ defined as above, $S$ is a Gr\"{o}bner-Shirshov
basis in $k \langle X;\mathcal{D}\rangle$ since in $S$ there are no
compositions except for the ambiguity $x_{i}x_{j}x_{k}$ which is a
trivial case. This implies that $A$ can be embedded into
$\mathcal{A}$. By (\ref{17})-(\ref{20}), $\mathcal{A}$ is a simple
associative differential algebra generated by $\{a,b\}$.\hfill
$\blacksquare$

\section{Associative algebras with multiple
operations}

Composition-Diamond lemma for associative algebra with multiple
operations $\Omega$ (associative $\Omega$-algebra, for short) is
established in a recent paper \cite{BoCQ}. By applying this lemma,
we show in this section that: (i). Every countably generated
associative $\Omega$-algebra can be embedded into a two-generated
associative $\Omega$-algebra. (ii). Any associative $\Omega$-algebra
can be embedded into a simple  associative $\Omega$-algebras. (iii).
Each countably generated  associative $\Omega$-algebra with
countable multiple operations $\Omega$ over a countable field $k$
can be embedded into a simple two-generated associative
$\Omega$-algebra.

The concept of multi-operations algebras ($\Omega$-algebras) was
first introduced by A.G. Kurosh  in \cite{Ku, Ku2}.

Let $k$ be a field. An associative algebra with multiple linear
operations is an associative $k$-algebra $A$ with a set $\Omega$ of
multi-linear operations.

Let $X$ be a set and
$$
\Omega=\bigcup_{n=1}^{\infty}\Omega_{n},
 $$
where $\Omega_{n}$ is the set of $n$-ary operations, for example,
ary $(\delta)=n$ if $\delta\in \Omega_n$.

Denote by   $S(X)$ the free semigroup without identity generated by
$X$.

For any non-empty set $Y$ (not necessarily a subset of $S(X)$), let
$$
\Omega(Y)=\bigcup\limits_{n=1}^{\infty}\{\delta(x_1,x_2,\cdots,x_n)|\delta\in
\Omega_n, x_i\in Y, \ i=1,2,\cdots,n\}.
$$

Define
\begin{eqnarray*}
\mathfrak{S}_{0}&=&S(X),\\
\mathfrak{S}_{1}&=&S(X\cup \Omega(\mathfrak{S}_{0})), \\
\vdots\ \ & &\ \ \ \ \vdots\\
\mathfrak{S}_{n}&=&S(X\cup \Omega(\mathfrak{S}_{n-1})),\\
\vdots\ \ & &\ \ \ \ \vdots
\end{eqnarray*}
Then we have
$$
\mathfrak{S}_{0}\subset\mathfrak{S}_{1}\subset\cdots
\subset\mathfrak{S}_{n}\subset\cdots.
$$
Let
$$
\mathfrak{S}(X)=\bigcup_{n\geq0}\mathfrak{S}_{n}.
$$
Then, we can see that $\mathfrak{S}(X)$ is a semigroup such that $
\Omega(\mathfrak{S}(X))\subseteq \mathfrak{S}(X). $

 For any $u\in \mathfrak{S}(X)$, $dep(u)=\mbox{min}\{n|u\in\mathfrak{S}_{n} \}$
 is called the depth of
 $u$.

Let   $k\langle X; \Omega\rangle$ be the $k$-algebra spanned by
$\mathfrak{S}(X)$.  Then, the element in $\mathfrak{S}(X)$ (resp.
$k\langle X; \Omega\rangle$) is called a $\Omega$-word (resp.
$\Omega$-polynomial).

Extend linearly each map $\delta\in\Omega_n$,
$$\delta:\mathfrak{S}(X)^n\rightarrow \mathfrak{S}(X), \
(x_1,x_2,\cdots,x_n)\mapsto \delta(x_1,x_2,\cdots,x_n)
$$
to $k\langle X; \Omega\rangle$. Then, $k\langle X; \Omega\rangle$ is
a free associative algebra with multiple linear operators $\Omega$
on the set $X$ (see \cite{BoCQ}).

Let $X$ and $\Omega$ be well ordered sets. We order $X^*$ by the
deg-lex ordering. For any $u\in \mathfrak{S}(X)$, $u$ can be
uniquely expressed without brackets as
$$
u=u_0\delta_{i_{_{1}}}\overrightarrow{x_{i_1}}u_1\cdots
\delta_{i_{_{t}}}\overrightarrow{x_{i_{t}}}u_t,
$$
where each $u_i\in X^*,\delta_{i_{_{k}}}\in \Omega_{i_{_{k}}}, \
\overrightarrow{x_{i_k}}=( x_{k_{1}},x_{k_{2}},\cdots,
x_{k_{i_{_{k}}}} )\in \mathfrak{S}(X)^{i_k}$.

Denote by
$$
wt(u)=(t,\delta_{i_{_{1}}},\overrightarrow{x_{i_{1}}},
\cdots,\delta_{i_{_{t}}},\overrightarrow{x_{i_t}}, u_0, u_1, \cdots,
u_t ).
$$
Then, we order  $\mathfrak{S}(X)$ as follows: for any $u,v\in
\mathfrak{S}(X)$,
\begin{equation}\label{o1}
u>v\Longleftrightarrow wt(u)>wt(v)\ \mbox{ lexicographically}
\end{equation}
by induction on $dep(u)+dep(v)$.

It is clear that the ordering (\ref{o1}) is a monomial ordering on
$\mathfrak{S}(X)$ (see \cite{BoCQ}).

Denote by $deg_{\Omega}(u)$ the number of $\delta$ in $u$ where
$\delta\in\Omega$, for example, if
$u=x_1\delta_1(x_2)\delta_3(x_2,x_1,\delta_1(x_3))$, then
$deg_{\Omega}(u)=3$.

\begin{theorem}
Every countably generated  associative $\Omega$-algebra can be
embedded into a two-generated associative $\Omega$-algebra.
\end{theorem}

\textbf{Proof}\   Suppose that $A= k\langle X ; \Omega |S\rangle$ is
an associative $\Omega$-algebra generated by $X$ with relations $S$,
where $X=\{x_{i},i=1,2,\dots\}$. By Shirshov algorithm, we can
assume that $S$ is a Gr\"{o}bner-Shirshov basis of the free
associative $\Omega$-algebra $k\langle X ; \Omega\rangle$ in the
sense of the paper \cite{BoCQ} with the ordering (\ref{o1}). Let
$H=k\langle X, a,b ;\Omega | S, aab^{i}ab=x_{i}, \
i=1,2,\dots\rangle$. We can check that $\{S, aab^{i}ab=x_{i}, \
i=1,2,\dots\}$ is a Gr\"{o}bner-Shirshov basis in the free
associative $\Omega$-algebra $k\langle X,a,b ; \Omega \rangle$ since
there are no new compositions. By the Composition-Diamond lemma in
\cite{BoCQ}, $A$ can be embedded into $H$ which is generated by
$\{a,b\}$.  \hfill $\blacksquare$

\begin{theorem}
Every  associative $\Omega$-algebra can be embedded into a simple
associative $\Omega$-algebra.
\end{theorem}

\textbf{Proof}\ Let  $A$ be an associative $\Omega$-algebra over a
field $k$ with $k$-basis $X=\{x_i\mid i\in I\}$ where $I$ is a well
ordered set. Denote by
\begin{eqnarray*}
S&=&\{x_ix_j=\{x_i,x_j\}, \delta_n(x_{k_1},\ldots,x_{k_n})=
\{\delta_n(x_{k_1},\ldots,x_{k_n})\}|\\
&& \ \ \ \ \ \ \ \ \ i,j,k_1,\ldots,k_n\in I,\
\delta_n\in\Omega_n,n\in N\},
\end{eqnarray*}
where $\{\delta_n(x_{k_1},\ldots,x_{k_n})\}$ is a linear combination
of $x_i, i\in I$. Then in the sense of the paper \cite{BoCQ}, $S$ is
a Gr\"{o}bner-Shirshov basis in the free associative
$\Omega$-algebra $k\langle X;\Omega \rangle$ with the ordering
(\ref{o1}). Therefore $A$ can be expressed as
\begin{eqnarray*}
A=k\langle X;\Omega |S\rangle.
\end{eqnarray*}
Let us totally order the set of monic elements of $A$. Denote by $T$
the set of indices for the resulting totally ordered set. Consider
the totally ordered set $T^2=\{(\theta,\sigma)\}$ and assign
$(\theta,\sigma)<(\theta',\sigma')$ if either $\theta<\theta'$ or
$\theta=\theta'$ and $\sigma<\sigma'$. Then $T^2$ is also totally
ordered set.

For each ordered pair of elements $f_\theta, f_\sigma\in\ A, \
\theta, \sigma\in T$, introduce the letters
$x_{\theta\sigma},y_{\theta\sigma}$.

Let $A_{1}$ be the associative $\Omega$-algebra given by the
generators
$$
X_1=\{x_i,y_{\theta\sigma},x_{\varrho\tau}\mid i\in I, \ \theta,
\sigma, \varrho, \tau\in T\}
$$
and the defining relations $S_1$ where $S_1$ is the union of $S$ and
\begin{eqnarray*}
x_{\theta\sigma}f_\theta y_{\theta\sigma}=f_\sigma,\
 ({\theta,\sigma})\in T^2.
\end{eqnarray*}
We can have that in $S_1$ there are no compositions unless for the
ambiguity $x_{i}x_{j}x_{k}$. But this case is trivial. Hence $S_1$
is a Gr\"{o}bner-Shirshov basis of the free associative
$\Omega$-algebra $k\langle X_1; \Omega \rangle$ in the sense of the
paper \cite{BoCQ} with the ordering (\ref{o1}). Thus, by the
Composition-Diamond lemma in \cite{BoCQ}, $A$ can be embedded into
$A_1$. The relations $x_{\theta\sigma}f_\theta
y_{\theta\sigma}=f_\sigma$ of $A_1$ provide that in ${A}_1$ every
monic element $f_\theta$ of the subalgebra ${A}$ generates an ideal
containing algebra $A$.

Mimicking the construction of the associative $\Omega$-algebra $A_1$
from the $A$, produce the associative $\Omega$-algebra $A_2$ from
$A_1$ and so on. As a result, we acquire an ascending chain of
associative $\Omega$-algebras $ A=A_0\subset A_1\subset
A_2\subset\cdots$ such that every  nonzero element generates the
same ideal. Let $ {\cal A}=\bigcup_ {k=0}^ {\infty}A_k. $ Then
${\cal A}$ is a simple associative $\Omega$-algebra. \hfill
$\blacksquare$

\begin{theorem}
Every countably generated associative $\Omega$-algebra with
countable multiple operations $\Omega$ over a countable field $k$
can be embedded into a simple two-generated associative
$\Omega$-algebra.
\end{theorem}
\textbf{Proof} \ Let $A$ be a countably generated associative
$\Omega$-algebra with countable multiple operations $\Omega$ over a
countable field $k$. We may assume that $A$ has a countable
$k$-basis $X_0=\{x_{i}| i=1,2,\ldots\}$. Denote by
\begin{eqnarray*}
S&=&\{x_ix_j=\{x_i,x_j\}, \delta_n(x_{k_1},\ldots,x_{k_n})=
\{\delta_n(x_{k_1},\ldots,x_{k_n})\}|\\
&& \ \ \ \ \ \ \ \ \ i,j,k_1,\ldots,k_n\in N,\
\delta_n\in\Omega_n,n\in N\}.
\end{eqnarray*}
Then  $A$ can be expressed as $ A=k\langle X_0;\Omega |S\rangle. $

Let $A_0=k\langle X_0;\Omega\rangle$, $A_0^+=A_0\backslash\{0\}$ and
fix the bijection
$$
(A_0^+,A_0^+)\longleftrightarrow\{(x_m^{(1)},y_m^{(1)}), m\in N\}.
$$

Let $X_1=X_0\cup\{x_m^{(1)},y_m^{(1)},a,b|m\in N\}$, $A_1=k\langle
X_1;\Omega\rangle$, $A_1^+=A_1\backslash\{0\}$ and fix the bijection
$$
(A_1^+,A_1^+)\longleftrightarrow\{(x_m^{(2)},y_m^{(2)}), m\in N\}.
$$
$$
\vdots
$$
Let $X_{n+1}=X_n\cup\{x_m^{(n+1)},y_m^{(n+1)}|m\in N\}$, $n\geq1$,
$A_{n+1}=k\langle X_{n+1};\Omega\rangle$,
$A_{n+1}^+=A_{n+1}\backslash\{0\}$ and fix the bijection
$$
(A_{n+1}^+,A_{n+1}^+)\longleftrightarrow\{(x_m^{(n+2)},y_m^{(n+2)}),
m\in N\}.
$$
$$
\vdots
$$

Consider the chain of the free associative $\Omega$-algebras
$$
A_0\subset A_1\subset A_2\subset\ldots\subset A_n\subset\ldots.
$$

Let $X=\bigcup_{n=0}^{\infty}X_n$. Then $k\langle
X;\Omega\rangle=\bigcup_ {n=0}^ {\infty}A_n$.

Now, define the desired algebra $\mathcal{A}$. Take the set $X$ as
the set of the generators for this algebra and take the union of $S$
and the following relations as one part of the relations for this
algebra
\begin{equation}\label{a2}
aa(ab)^nb^{2m+1}ab=x_m^{(n)},\  m,n\in N
\end{equation}
\begin{equation}\label{a3}
aa(ab)^nb^{2m}ab=y_m^{(n)}, \ m,n\in N
\end{equation}
\begin{equation}\label{a4}
aabbab=x_1
\end{equation}

Before we introduce the another part of the relations on
$\mathcal{A}$, let us define canonical words of the algebras $A_n$,
$n\geq0$. A  $\Omega$-word in $X_0$ without subwords that are the
leading terms of $s\ (s\in S)$ is called a canonical word of $A_0$.
A $\Omega$-word in $X_1$ without subwords that are the leading terms
of  $s\ (s\in S\cup\{(\ref{a2}), (\ref{a3}), (\ref{a4})\})$ and
without subwords of the form
$$
(\delta_1(x_m^{(1)}))^{deg_{\Omega}(\overline{g^{(0)}})+1}\overline{f^{(0)}}y_m^{(1)},
$$
where
$(x_m^{(1)},y_m^{(1)})\longleftrightarrow(f^{(0)},g^{(0)})\in(A_0^+,A_0^+)$
such that $f^{(0)},g^{(0)}$ are non-zero linear combination of
canonical words of $A_0$, is called a canonical word of $A_1$.
Suppose that we have defined canonical word of $A_{k}$, $k<n$. A
$\Omega$-word in $X_n$ without subwords that are the leading terms
of $s\ (s\in S\cup\{(\ref{a2}), (\ref{a3}), (\ref{a4})\})$ and
without subwords of the form
$$
(\delta_1(x_m^{(k+1)}))^{deg_{\Omega}(\overline{g^{(k)}})+1}\overline{f^{(k)}}y_m^{(k+1)},
$$
where
$(x_m^{(k+1)},y_m^{(k+1)})\longleftrightarrow(f^{(k)},g^{(k)})\in(A_k^+,A_k^+)$
such that $f^{(k)},g^{(k)}$ are non-zero linear combination of
canonical words of $A_k$, is called a canonical word of $A_n$.

Then the another part of the relations on $\mathcal{A}$ are the
following:
\begin{equation}\label{a5}
(\delta_1(x_m^{(n)}))^{deg_{\Omega}(\overline{g^{(n-1)}})+1}f^{(n-1)}y_m^{(n)}-g^{(n-1)}=0,
\ \ m, n\in N
\end{equation}
where
$(x_m^{(n)},y_m^{(n)})\longleftrightarrow(f^{(n-1)},g^{(n-1)})\in(A_{n-1}^+,A_{n-1}^+)$
such that $f^{(n-1)},g^{(n-1)}$ are non-zero linear combination of
canonical words of $A_{n-1}$.

We can see that in $\mathcal{A}$ every element can be expressed as
linear combination of canonical words.

Denote by $S_1=S\cup\{(\ref{a2}), (\ref{a3}),
(\ref{a4}),(\ref{a5})\}$. We can have that with the ordering
(\ref{o1}), $S_1$ is a Gr\"{o}bner-Shirshov basis in $k \langle
X;\Omega\rangle$ in the sense of the paper \cite{BoCQ} since in
$S_1$ there are no compositions except for the ambiguity
$x_{i}x_{j}x_{k}$ which is a trivial case. This implies that $A$ can
be embedded into $\mathcal{A}$. By (\ref{a2})-(\ref{a5}),
$\mathcal{A}$ is a simple associative $\Omega$-algebra generated by
$\{a,b\}$. \hfill $\blacksquare$

\section{Associative $\lambda$-differential algebras}

In this section, by applying the Composition-Diamond lemma for
associative $\Omega$-algebras in \cite{BoCQ}, we show that: (i).
Each countably generated associative $\lambda$-differential algebra
can be embedded into a two-generated associative
$\lambda$-differential algebra. (ii). Each associative
$\lambda$-differential algebra can be embedded into a simple
associative $\lambda$-differential algebra. (iii). Each countably
generated associative $\lambda$-differential algebra over a
countable field $k$ can be embedded into a simple two-generated
associative $\lambda$-differential algebra.

Let $k$ be a commutative ring with unit and $\lambda\in k$. An
associative $\lambda$-differential algebra over $k$ (\cite{GuK08})
is an  associative $k$-algebra  $R$ together with a  $k$-linear
operator $D:R\rightarrow R$ such that
$$
D(xy)=D(x)y+xD(y)+\lambda D(x)D(y),\ \forall x, y \in R.
$$

Any  associative $\lambda$-differential  algebra is also an
associative algebra with one operator $\Omega=\{D\}$.

In this section, we will use the notations given in the Section 6.

Let $X$ be well ordered and $k\langle X;D\rangle$ the free
associative algebra with one operator  $\Omega=\{D\}$ defined in the
Section 6.

 For any $u\in
\mathfrak{S}(X)$, $u$ has a unique expression
$$
u=u_1u_2\cdots u_n,
$$
where each $u_i\in X\cup D(\mathfrak{S}(X))$. Denote by
$deg_{_{X}}(u)$ the number of $x\in X$ in $u$, for example, if
$u=D(x_1x_2)D(D(x_1))x_3\in \mathfrak{S}(X)$, then
$deg_{_{X}}(u)=4$.  Let
$$
wt(u)=(deg_{_{X}}(u), u_1,u_2,\cdots, u_n).
$$
Now, we order $\mathfrak{S}(X)$ as follows: for any $u,v\in
\mathfrak{S}(X)$,
\begin{equation}\label{o2}
u>v\Longleftrightarrow wt(u)>wt(v)\ \mbox{ lexicographically}
\end{equation}
where for each  $t, \ u_t>v_t$ if one of the following holds:

(a) $u_t, v_t\in X$ and $u_t>v_t$;

(b)  $u_t=D(u_{t}^{'}), v_t\in X$;

(c)  $u_t=D(u_{t}^{'}),v_t=D(v_{t}^{'})$ and $u_{t}^{'}>v_{t}^{'}$.

\ \

Then the ordering (\ref{o2}) is a monomial ordering on
$\mathfrak{S}(X)$ (see \cite{BoCQ}).

\begin{lemma}(\cite{BoCQ}, Theorem 5.1)\label{t3.5}
With  the  ordering (\ref{o2})  on $\mathfrak{S}(X)$,
$$
S_0=\{D(uv)-D(u)v-uD(v)-\lambda D(u)D(v) |\  u,v \in
\mathfrak{S}(X)\}
$$
is a Gr\"{o}bner-Shirshov basis in the free $\Omega$-algebra
$k\langle X;D\rangle$ where $\Omega=\{D\}$.
\end{lemma}

\begin{lemma}\label{l4.4}
Let $A$ be an associative $\lambda$-differential algebra with
$k$-basis $X=\{x_{i}| i\in I\}$. Then $A$ has a representation $A=
k\langle X ; D | S \rangle$, where
$S=\{x_{i}x_{j}=\{x_{i},x_{j}\},D(x_{i})=\{D(x_{i})\},
D(x_{i}x_{j})=D(x_{i})x_{j}+x_{i}D(x_{j})+\lambda
D(x_{i})D(x_{j})\mid  i, j\in I\}$.

Moreover, if $I$ is a well ordered set, then with the ordering
(\ref{o2}) on $\mathfrak{S}(X)$, $S$ is a Gr\"{o}bner-Shirshov basis
in the free $\Omega$-algebra $k\langle X; D\rangle$  in the sense of
\cite{BoCQ}.
\end{lemma}
\textbf{Proof}\ Clearly, $k\langle X ; D | S \rangle$ is an
associative $\lambda$-differential  algebra.  By the
Composition-Diamond lemma in \cite{BoCQ}, it suffices to check that
with the ordering (\ref{o2}) on $\mathfrak{S}(X)$,  $S$ is a
Gr\"{o}bner-Shirshov basis in $k\langle X; D\rangle$  in the sense
of \cite{BoCQ}.

The ambiguities $w$ of all possible compositions of
$\Omega$-polynomials in $S$ are:
\begin{enumerate}
\item[1)]\ $x_ix_jx_k,\ i,j,k\in I,$
\item[2)]\ $D(x_ix_j),\ i,j\in I$.
\end{enumerate}
We will check that each composition in $S$ is trivial $mod(S,w)$.

For 1), the result is trivial.

For 2), let $f=D(x_{i}x_{j})-D(x_{i})x_{j}-x_{i}D(x_{j})-\lambda
D(x_{i})D(x_{j})$,\ \ $g=x_{i}x_{j}-\{x_{i},x_{j}\},\ i,j\in I$.
Then $w=D(x_{i}x_{j})$ and
\begin{eqnarray*}
(f,g)_{w}&=&-D(x_{i})x_{j}-x_{i}D(x_{j})-\lambda
D(x_{i})D(x_{j})+D(\{x_{i},x_{j}\})\\
&\equiv&\{D(\{x_{i},x_{j}\})\}-\{\{D(x_{i})\},x_{j}\}-\{x_{i},\{D(x_{j})\}\}-\lambda
\{\{D(x_{i})\},\{D(x_{j})\}\}\\
&\equiv&0 \ \ mod(S,w).
\end{eqnarray*}
This shows that $S$ is a Gr\"{o}bner-Shirshov basis in the free
$\Omega$-algebra $k\langle X; D\rangle$. \hfill $\blacksquare$

 \ \

Now we get the embedding theorems for associative
$\lambda$-differential algebras.

\begin{theorem} \label{t3.6}
Every countably generated associative $\lambda$-differential algebra
over a field can be embedded into a two-generated associative
$\lambda$-differential algebra.
\end{theorem}

\textbf{Proof}\ Let $A$ be a  countably generated associative
$\lambda$-differential algebra over a field $k$. We may assume that
$A$ has a countable $k$-basis $X=\{x_{i}| i=1,2,\ldots\}$. By Lemma
\ref{l4.4}, $A= k\langle X ; D | S \rangle$, where
$S=\{x_{i}x_{j}=\{x_{i},x_{j}\},D(x_{i})=\{D(x_{i})\},
D(x_{i}x_{j})=D(x_{i})x_{j}+x_{i}D(x_{j})+\lambda
D(x_{i})D(x_{j})\mid  i, j\in N\}$.

Let $H=k\langle X, a,b; D | S_{1} \rangle$ where
\begin{eqnarray*}
&&S_{1}=\{x_{i}x_{j}=\{x_{i},x_{j}\},D(x_{i})=\{D(x_{i})\},
D(uv)=D(u)v+uD(v)+\lambda D(u)D(v),\\
&&\ \ \ \ \ \ \ \ aab^{i}ab=x_{i}| u, v \in \mathfrak{S}(X,a,b), i,
j\in N\}.
\end{eqnarray*}

We want to prove that $S_{1}$ is also a Gr\"{o}bner-Shirshov basis
in the free $\Omega$-algebra $k\langle X,a,b;D\rangle$ with the
ordering (\ref{o2}). Now, let us check all the possible compositions
in $S_{1}$. The ambiguities $w$ of  all possible compositions of
$\Omega$-polynomials in $S_{1}$ are: \\
$$
\begin{array}{lllll}
1)\ \ x_{i}x_{j}x_{k}& 2)\ D(u|_{x_{i}x_{j}}v)& 3)\
D(uv|_{x_{i}x_{j}})& 4)\ D(u|_{D(x_{i})}v)& 5)\ D(uv|_{D(x_{i})})\\
6)\ D(uv|_{D(u_{1}v_{1})})& 7) \ D(u|_{D(u_{1}v_{1})}v)& 8)\
D(u|_{aab^{i}ab}v)& 9)\ D(uv|_{aab^{i}ab})
\end{array}
$$
where $u,v,u_{1},v_{1}\in \mathfrak{S}(X,a,b),x_{i},x_{j},x_{k}\in X
$. We have to check that all these compositions are trivial
$mod(S_1,w)$. In fact, by Lemma \ref{t3.5} and since $S$ is a
Gr\"{o}bner-Shirshov basis in $k\langle X ;D\rangle$, we need only
to check $2)-5),8),9)$. Here, for example, we just check $3),4),8)$.
Others are similarly proved.

For 3), let
$f=D(uv|_{x_{i}x_{j}})-D(u)v|_{x_{i}x_{j}}-uD(v|_{x_{i}x_{j}})-\lambda
D(u)D(v|_{x_{i}x_{j}}), \ \ g=x_{i}x_{j}-\{x_{i},x_{j}\}, \ \ u,
v\in \mathfrak{S}(X,a,b), x_{i},x_{j}\in X  $. Then
$w=D(uv|_{x_{i}x_{j}})$ and
\begin{eqnarray*}
(f,g)_{w}&=&- D(u)v|_{x_{i}x_{j}} - uD(v|_{x_{i}x_{j}})-
\lambda D(u)D(v|_{x_{i}x_{j}})+D(uv|_{\{x_{i},x_{j}\}})\\
&\equiv&- D(u)v|_{\{x_{i},x_{j}\}} - uD(v|_{\{x_{i},x_{j}\}})-
\lambda D(u)D(v|_{\{x_{i},x_{j}\}})+D(uv|_{\{x_{i},x_{j}\}})\\
&\equiv&0.
\end{eqnarray*}

For 4), let
$f=D(u|_{D(x_{i})}v)-D(u|_{D(x_{i})})v-u|_{D(x_{i})}D(v)- \lambda
D(u|_{D(x_{i})})D(v), \ \ g=D(x_{i})-\{D(x_{i})\},\ \
u,v,D(x_{i})\in \mathfrak{S}(X,a,b),x_{i}\in X $. Then
$w=D(u|_{D(x_{i})}v)$ and
\begin{eqnarray*}
(f,g)_{w}&=&- D( u|_{D(x_{i})})v - u|_{D(x_{i})}D(v)-
\lambda D(u|_{D(x_{i})})D(v)+D(u|_{\{D(x_{i})\}}v)\\
&\equiv&- D( u|_{\{D(x_{i})\}})v - u|_{\{D(x_{i})\}}D(v)-
\lambda D(u|_{\{D(x_{i})\}})D(v)+D(u|_{\{D(x_{i})\}}v)\\
&\equiv&0.
\end{eqnarray*}

For 8), let $f=D(u|_{aab^{i}ab}v) - D( u|_{aab^{i}ab})v -
u|_{aab^{i}ab}D(v)-\lambda D(u|_{aab^{i}ab})D(v) , \
g=aab^{i}ab-x_{i},\ u, v \in \mathfrak{S}(X,a,b),\ x_{i}\in X$. Then
$w=D(u|_{aab^{i}ab}v)$ and
\begin{eqnarray*}
(f,g)_{w}&=&- D( u|_{aab^{i}ab})v - u|_{aab^{i}ab}D(v)-
\lambda D(u|_{aab^{i}ab})D(v)+D(u|_{x_{i}}v)\\
&\equiv&D(u|_{x_{i}}v)- D( u|_{x_{i}})v - u|_{x_{i}}D(v)- \lambda
D(u|_{x_{i}})D(v)\\
&\equiv&0.
\end{eqnarray*}

So  $S_{1}$ is a Gr\"{o}bner-Shirshov basis in $k\langle
X,a,b;D\rangle$.  By the Composition-Diamond lemma in \cite{BoCQ},
$A$ can be embedded into $H$ which is generated by $\{a,b\}$. \hfill
$\blacksquare$

\begin{theorem}\label{t3.7}
Every associative $\lambda$-differential algebra over a field can be
embedded into a simple associative $\lambda$-differential algebra.
\end{theorem}

\textbf{Proof}\ Let $A$ be an associative $\lambda$-differential
algebra over a field $k$ with basis $X=\{x_i\mid i\in I\}$ where $I$
is a well ordered set. Then by Lemma \ref{l4.4}, $A$ can be
expressed as $A= k\langle X ; D | S \rangle$ where
$S=\{x_{i}x_{j}=\{x_{i},x_{j}\},D(x_{i})=\{D(x_{i})\},
D(x_{i}x_{j})=D(x_{i})x_{j}+x_{i}D(x_{j})+\lambda
D(x_{i})D(x_{j})\mid  i, j\in I\}$ and $S$ is a Gr\"{o}bner-Shirshov
basis in $k \langle X ;D\rangle$  with the ordering (\ref{o2}). Let
us totally order the set of monic elements of ${A}$. Denote by $T$
the set of indices for the resulting totally ordered set. Consider
the totally ordered set $T^2=\{(\theta,\sigma)\}$ and assign
$(\theta,\sigma)<(\theta',\sigma')$ if either $\theta<\theta'$ or
$\theta=\theta'$ and $\sigma<\sigma'$. Then $T^2$ is also totally
ordered set.

For each ordered pair of elements $f_\theta, f_\sigma\in\ A,\
\theta, \sigma\in T$, introduce the letters
$x_{\theta\sigma},y_{\theta\sigma}$.

Let $A_{1}$ be the associative $\lambda$-differential algebra given
by the generators
$$
X_1=\{x_i, y_{\theta\sigma}, x_{\varrho\tau} | i\in I, \ \theta,
\sigma, \varrho, \tau\in T\}
$$
and the defining relations
$$
x_ix_j=\{x_i,x_j\}, \ \ i,j\in I,
$$
$$
D(x_{i})=\{D(x_{i})\}, \ \ i\in I,
$$
$$
D(uv)=D(u)v+uD(v)+\lambda D(u)D(v),\ \  u,v \in \mathfrak{S}(X_1),
$$
$$
x_{\theta\sigma}f_\theta y_{\theta\sigma}=f_\sigma, \ \
({\theta,\sigma})\in T^2.
$$
 We want to prove
that these relations is also a Gr\"{o}bner-Shirshov basis in $k
\langle X_1 ;D\rangle$ with the same ordering (\ref{o2}). Now, let
us check all the possible compositions. The ambiguities $w$ of  all
possible compositions of
$\Omega$-polynomials are: \\
$$
\begin{array}{lllll}
1) \  x_{i}x_{j}x_{k}& 2)\ D(u|_{x_{i}x_{j}}v)& 3)\
D(uv|_{x_{i}x_{j}})& 4)\ D(u|_{D(x_{i})}v)&
5)\ D(uv|_{D(x_{i})})\\
6)\ D(uv|_{D(u_{1}v_{1})})& 7)\  D(u|_{D(u_{1}v_{1})}v)& 8)\
D(u|_{x_{\theta\sigma}\overline{f_\theta} y_{\theta\sigma}}v)& 9)\
D(uv|_{x_{\theta\sigma}\overline{f_\theta} y_{\theta\sigma}})
\end{array}
$$
where $u,v,u_{1},v_{1}\in \mathfrak{S}(X_1), x_{i}, x_{j}, x_{k}\in
X, ({\theta,\sigma})\in T^2$.

In fact, by  Lemma \ref{t3.5} and since $S$ is a
Gr\"{o}bner-Shirshov basis in $K\langle X ;D\rangle$, we just need
to check $2)-5),8),9)$. Here, for example, we just check $8)$.
Others are similarly proved.

Let $f=D(u|_{x_{\theta\sigma}\overline{f_\theta }y_{\theta\sigma}}v)
- D( u|_{x_{\theta\sigma}\overline{f_\theta} y_{\theta\sigma}})v -
u|_{x_{\theta\sigma}\overline{f_\theta}
y_{\theta\sigma}}D(v)-\lambda
D(u|_{x_{\theta\sigma}\overline{f_\theta} y_{\theta\sigma}})D(v), \
g=x_{\theta\sigma}f_\theta
y_{\theta\sigma}-f_\sigma=x_{\theta\sigma}\overline{f_\theta
}y_{\theta\sigma}+x_{\theta\sigma}f_\theta'y_{\theta\sigma}-f_\sigma$,
where $f_\theta=\overline{f_\theta }+f_\theta'$, $u, v \in
\mathfrak{S}(X_1), ({\theta,\sigma})\in T^2$. Then
$w=D(u|_{x_{\theta\sigma}\overline{f_\theta }y_{\theta\sigma}}v)$
and
\begin{eqnarray*}
(f,g)_{w}&=&- D( u|_{x_{\theta\sigma}\overline{f_\theta}
y_{\theta\sigma}})v - u|_{x_{\theta\sigma}\overline{f_\theta}
y_{\theta\sigma}}D(v)-\lambda
D(u|_{x_{\theta\sigma}\overline{f_\theta} y_{\theta\sigma}})D(v)+
D(u|_{(-x_{\theta\sigma}f_\theta'y_{\theta\sigma}+f_\sigma)}v)\\
&\equiv&D(u|_{(-x_{\theta\sigma}f_\theta'y_{\theta\sigma}+f_\sigma)}v)-
D( u|_{(-x_{\theta\sigma}f_\theta'y_{\theta\sigma}+f_\sigma)})v -
u|_{(-x_{\theta\sigma}f_\theta'y_{\theta\sigma}+f_\sigma)}D(v)\\
&&- \lambda
D(u|_{(-x_{\theta\sigma}f_\theta'y_{\theta\sigma}+f_\sigma)})D(v)\\
&\equiv&0.
\end{eqnarray*}
Thus, by the Composition-Diamond lemma in \cite{BoCQ}, $A$ can be
embedded into $A_1$. The relations $x_{\theta\sigma}f_\theta
y_{\theta\sigma}=f_\sigma$ of $A_1$ provide that in ${A}_1$ every
monic element $f_\theta$ of the subalgebra ${A}$ generates an ideal
containing algebra $A$.

Mimicking the construction of the associative $\lambda$-differential
algebra $A_1$ from the $A$, produce the associative
$\lambda$-differential algebra $A_2$ from $A_1$ and so on. As a
result, we acquire an ascending chain of associative
$\lambda$-differential algebras $ A=A_0\subset A_1\subset
A_2\subset\dots$ such that every  nonzero element generates the same
ideal. Let $ {\cal A}=\bigcup_ {k=0}^ {\infty}A_k. $ Then ${\cal A}$
is a simple associative $\lambda$-differential algebra. \hfill
$\blacksquare$

\begin{theorem}
Every countably generated associative $\lambda$-differential algebra
over a countable field $k$ can be embedded into a simple
two-generated associative $\lambda$-differential algebra.
\end{theorem}
\textbf{Proof} \ Let $A$ be a countably generated associative
$\lambda$-differential algebra over a countable field $k$. We may
assume that $A$ has a countable $k$-basis $X_0=\{x_{i}|
i=1,2,\ldots\}$ and it can be expressed as, by Lemma \ref{l4.4}, $A=
k\langle X_0 ; D | S_0 \rangle$ where
$S_0=\{x_{i}x_{j}=\{x_{i},x_{j}\},D(x_{i})=\{D(x_{i})\},
D(x_{i}x_{j})=D(x_{i})x_{j}+x_{i}D(x_{j})+\lambda
D(x_{i})D(x_{j})\mid  i, j\in N\}$ and  $S_0$ is a
Gr\"{o}bner-Shirshov basis in $k \langle X_0 ;D\rangle$ with the
ordering (\ref{o2}).

Let $A_0=k\langle X_0;D \rangle$, $A_0^+=A_0\backslash\{0\}$ and fix
the bijection
$$
(A_0^+,A_0^+)\longleftrightarrow\{(x_m^{(1)},y_m^{(1)}), m\in N\}.
$$

Let $X_1=X_0\cup\{x_m^{(1)},y_m^{(1)},a,b|m\in N\}$, $A_1=k\langle
X_1;D\rangle$, $A_1^+=A_1\backslash\{0\}$ and fix the bijection
$$
(A_1^+,A_1^+)\longleftrightarrow\{(x_m^{(2)},y_m^{(2)}), m\in N\}.
$$
$$
\vdots
$$
Let $X_{n+1}=X_n\cup\{x_m^{(n+1)},y_m^{(n+1)}|m\in N\}$, $n\geq1$,
$A_{n+1}=k\langle X_{n+1};D \rangle$,
$A_{n+1}^+=A_{n+1}\backslash\{0\}$ and fix the bijection
$$
(A_{n+1}^+,A_{n+1}^+)\longleftrightarrow\{(x_m^{(n+2)},y_m^{(n+2)}),
m\in N\}.
$$
$$
\vdots
$$

Consider the chain of the free $\Omega$-algebras
$$
A_0\subset A_1\subset A_2\subset\ldots\subset A_n\subset\ldots.
$$

Let $X=\bigcup_ {n=0}^ {\infty}X_n$. Then $k\langle
X;D\rangle=\bigcup_ {n=0}^ {\infty}A_n$.

Now, define the desired algebra $\mathcal{A}$. Take the set $X$ as
the set of the generators for this algebra and take the following
relations as one part of the relations for this algebra
\begin{equation}\label{b11}
x_{i}x_{j}=\{x_{i},x_{j}\}, D(x_{i})=\{D(x_{i})\},\ i,j\in N
\end{equation}
\begin{equation}\label{b12}
D(uv)=D(u)v+uD(v)+\lambda D(u)D(v),\ u,v \in \mathfrak{S}(X)
\end{equation}
\begin{equation}\label{b2}
aa(ab)^nb^{2m+1}ab=x_m^{(n)},\  m,n\in N
\end{equation}
\begin{equation}\label{b3}
aa(ab)^nb^{2m}ab=y_m^{(n)}, \ m,n\in N
\end{equation}
\begin{equation}\label{b4}
aabbab=x_1
\end{equation}

Before we introduce the another part of the relations on
$\mathcal{A}$, let us define canonical words of the algebras $A_n$,
$n\geq0$. A $\Omega$-word in $X_0$ without subwords that are the
leading terms of (\ref{b11}) and (\ref{b12}) is called a canonical
word of $A_0$. A $\Omega$-word in $X_1$ without subwords that are
the leading terms of (\ref{b11}), (\ref{b12}), (\ref{b2}),
(\ref{b3}), (\ref{b4}) and without subwords of the form
$$
(x_m^{(1)})^{deg_X(\overline{g^{(0)}})}\overline{f^{(0)}}y_m^{(1)},
$$
where
$(x_m^{(1)},y_m^{(1)})\longleftrightarrow(f^{(0)},g^{(0)})\in(A_0^+,A_0^+)$
such that $f^{(0)},g^{(0)}$ are non-zero linear combination of
canonical words of $A_0$, is called a canonical word of $A_1$.
Suppose that we have defined canonical word of $A_{k}$, $k<n$. A
$\Omega$-word in $X_n$ without subwords that are the leading terms
of (\ref{b11}), (\ref{b12}), (\ref{b2}), (\ref{b3}), (\ref{b4}) and
without subwords of the form
$$
(x_m^{(k+1)})^{deg_X(\overline{g^{(k)}})}\overline{f^{(k)}}y_m^{(k+1)},
$$
where
$(x_m^{(k+1)},y_m^{(k+1)})\longleftrightarrow(f^{(k)},g^{(k)})\in(A_k^+,A_k^+)$
such that $f^{(k)},g^{(k)}$ are non-zero linear combination of
canonical words of $A_k$, is called a canonical word of $A_n$.

Then the another part of the relations on $\mathcal{A}$ are the
following:
\begin{equation}\label{b5}
(x_m^{(n)})^{deg_X(\overline{g^{(n-1)}})}f^{(n-1)}y_m^{(n)}-g^{(n-1)}=0,\
\ m, n\in N
\end{equation}
where
$(x_m^{(n)},y_m^{(n)})\longleftrightarrow(f^{(n-1)},g^{(n-1)})\in(A_{n-1}^+,A_{n-1}^+)$
such that $f^{(n-1)},g^{(n-1)}$ are non-zero linear combination of
canonical words of $A_{n-1}$.

We can get that in $\mathcal{A}$ every element can be expressed as
linear combination of canonical words.

Denote by $S$ the set constituted by the relations
(\ref{b11})-(\ref{b5}). We want to prove that $S$ is also a
Gr\"{o}bner-Shirshov basis in the free $\Omega$-algebra $k\langle
X;D\rangle$ with the ordering (\ref{o2}).  The ambiguities $w$ of
all possible compositions of
$\Omega$-polynomials in $S$ are: \\
$$
\begin{array}{lll}
1)\ \ x_{i}x_{j}x_{k}& 2)\ D(u|_{x_{i}x_{j}}v)& 3)\
D(uv|_{x_{i}x_{j}})\\ 4)\ D(u|_{D(x_{i})}v)& 5)\ D(uv|_{D(x_{i})})&
6)\ D(uv|_{D(u_{1}v_{1})})\\ 7) \ D(u|_{D(u_{1}v_{1})}v)& 8)\
D(u|_{aa(ab)^nb^{2m+1}ab}v)& 9)\ D(uv|_{aa(ab)^nb^{2m+1}ab}) \\10)\
D(u|_{aa(ab)^nb^{2m}ab}v)& 11)\ D(uv|_{aa(ab)^nb^{2m}ab}) &12)\
D(u|_{aabbab}v)\\ 13)\ D(uv|_{aabbab})& 14)\
D(u|_{(x_m^{(n)})^{deg_X(\overline{g^{(n-1)}})}\overline{f^{(n-1)}}y_m^{(n)}}v)&
15)\
D(uv|_{(x_m^{(n)})^{deg_X(\overline{g^{(n-1)}})}\overline{f^{(n-1)}}y_m^{(n)}})
\end{array}
$$
where $u,v,u_{1},v_{1}\in \mathfrak{S}(X_1), x_{i}, x_{j}, x_{k}\in
X$. The proof of all possible compositions to be trivial $mod(S,w)$
is similar to that of Theorems \ref{t3.6}, \ref{t3.7}. Here we omit
the details. So  $S$ is a Gr\"{o}bner-Shirshov basis in $k\langle
X;D\rangle$ with the ordering (\ref{o2}), which implies that $A$ can
be embedded into $\mathcal{A}$. By (\ref{b2})-(\ref{b5}),
$\mathcal{A}$ is a simple associative $\lambda$-differential algebra
generated by $\{a,b\}$. \hfill $\blacksquare$

\section{ Modules}

In this section, by applying the Composition-Diamond lemma for
modules (see \cite{ChCZ,Ch}), we show that every countably generated
$k\langle X \rangle$-module can be embedded into a cyclic $k\langle
X \rangle$-module, where $|X|>1$.

Let $X,Y$ be well ordered sets and $mod_{k\langle X\rangle} \langle
Y\rangle$ a free left $k\langle X\rangle$-module with the basis $Y$.
Suppose that $<$ is the deg-lex ordering on $X^*$. Let
$X^*Y=\{uy|u\in X^*, \ y\in Y\}$. We define an ordering $\prec$ on
$X^*Y$ as follows: for any $w_1=u_1y_i,w_2=u_2y_j\in X^*Y$,
\begin{equation}\label{03}
w_1\prec w_2\Leftrightarrow u_1<u_2 \ \ \mbox{ or } u_1=u_2, \
y_i<y_j
\end{equation}
It is clear that the ordering $\prec$ is left
compatible in the sense of
$$
w\prec w'\Rightarrow aw \prec aw' \ \mbox{ for any  } a\in X^*.
$$

\begin{theorem}\label{t8.8}
Let $X$ be a set with $|X|>1$. Then every countably generated
$k\langle X \rangle$-module can be embedded into a cyclic $k\langle
X \rangle$-module.
\end{theorem}

\textbf{Proof}\  We may assume that  $M=Mod_{k\langle X
\rangle}\langle Y| T \rangle$ where $Y=\{y_{i},i=1,2,\dots\}$. By
Shirshov algorithm, we may assume that $T$ is a Gr\"{o}bner-Shirshov
basis in the free module $Mod_{k\langle X \rangle}\langle Y\rangle$
in the sense of the paper \cite{ChCZ} with the ordering (\ref{03})
on $X^*Y$.

Assume that $ a, b\in X,\ a\neq b$. Consider the ${k\langle X
\rangle}$-module
$$
_{{k\langle X \rangle}}M'=Mod_{k\langle X \rangle} \langle Y, y |
T,\ ab^{i}y-y_{i},y_{i}\in Y, i=1,2,\dots \rangle.
$$
 We can  check that $\{T,\  ab^{i}y-y_{i},i=1,2,\dots\}$ is
 also a Gr\"{o}bner-Shirshov basis
in the free module $Mod_{k\langle X \rangle} \langle Y, y  \rangle $
with the same ordering (\ref{03}) on $X^*(Y\cup\{y\})$ since there
are no new compositions. By the Composition-Diamond lemma in
\cite{ChCZ}, $M$ can be embedded into $_{{k\langle X \rangle}}M'$
which  is a cyclic $k\langle X \rangle$-module generated by $y$.
  \hfill
$\blacksquare$

\ \

\noindent\textbf{Remark}\ In Theorem \ref{t8.8}, the condition
$|X|>1$ is essential. For example, let $_{k[x]}M=\oplus_{i\in
I}k[x]y_i$ be a free $k[x]$-module with $k[x]$-basis $Y=\{y_i|i\in
I\}$, where $|I|>1$. Then $_{k[x]}M$ can not be embedded into a
cyclic $k[x]$-module. Indeed, suppose that $_{k[x]}M$ can be
embedded into a cyclic $k[x]$-module $k[x]y$. Let $y_1,y_2\in Y$
with $y_1\neq y_2$. Then there exist $f(x),g(x)\in k[x]$ such that
$y_1=f(x)y, \ y_2=g(x)y$. This implies that $g(x)y_1=f(x)y_2$, a
contradiction.


\begin{thebibliography}{9}

\bibitem{Al} A.A. Albert, A note on the exceptional Jordan algebra,
{\it Proc. Nat. Acad. Sci. U.S.A.}, {\bf 36}, 372-374 (1950).



\bibitem{Bax} G. Baxter, An analytic problem whose solution follows
from a simple algebraic identity, {\it Pacific J. Math.}, {\bf 10},
731-742 (1960).


\bibitem{Be} G.M. Bergman, The diamond lemma for ring theory,  {\it Adv. in
Math.}, {\bf29}, 178-218 (1978).

\bibitem{Bo62} L.A. Bokut, Embedding  algebras into algebraically
closed Lie algebras,  {\it Dokl. Akad. Nauk SSSR}, {\bf14}(5),
963-964 (1962).

\bibitem{Bo62L} L.A. Bokut, Embedding Lie algebras into algebraically
closed Lie algebras, {\it Algebra i Logika}, {\bf1}(2), 47-53
(1962).

\bibitem{Bo63} L.A. Bokut, Some embedding theorems for rings and semigroups. I,
II, {\it Sibirsk. Mat. Zh.,} {\bf4}, 500-518, 729-743 (1963).

\bibitem{Bo72} L.A. Bokut, Unsolvability of the word problem, and
subalgebras of finitely presented Lie algebras,  {\it Izv. Akad.
Nauk. SSSR Ser. Mat.},  {\bf36}, 1173-1219 (1972).

\bibitem{Bo76} L.A. Bokut, Imbeddings into simple associative
algebras,  {\it Algebra i Logika},  {\bf15}, 117-142 (1976).

\bibitem{Bo78}L.A. Bokut, Imbedding into algebraically closed and simple Lie algebras, {\it Trudy
Mat. Inst. Steklov.}, {\bf148}, 30-42 (1978).

\bibitem{BoC1} L.A. Bokut and Yuqun Chen, Gr\"{o}bner-Shirshov basis
for free Lie algebras: after A.I. Shirshov, {\it Southeast Asian
Bull. Math.}, {\bf31}, 1057-1076 (2007).  arXiv:0804.1254

\bibitem{BoCQ} L.A. Bokut, Yuqun Chen and Jianjun Qiu, Gr\"{o}bner-Shirshov basis
for Associative Algebras with Multiple Operator and Free Rota-Baxter
Algebras, {\it  Journal of Pure and Applied Algebra}, to appear.
doi:10.1016/j.jpaa.2009.05.005


\bibitem{BoK} L.A. Bokut and G.P. Kukin, Algorithmic and
Combinatorial Algebra, Kluwer Academic Publishers, 1994.

\bibitem{Bu65} B. Buchberger, An algorithm for finding a basis for the
residue class ring of a zero-dimensional polynomial ideal [in
German], {\it Ph.D. thesis}, University of Innsbruck, Austria
(1965).

\bibitem{Bu70} B. Buchberger, An algorithmical criteria for the
solvability of algebraic systems of equations [in German], {\it
Aequationes Math.}, {\bf4}, 374-383 (1970).


\bibitem{Chen} K.T. Chen, R.H. Fox and R.C. Lyndon, Free differential calculus, IV:
the quotient groups of the lower central series, {\it Annals of
Mathematics}, {\bf68}, 81-95 (1958).


\bibitem{ChCL} Yuqun Chen, Yongshan Chen and Yu Li, Composition-Diamond
Lemma for differential algebras, {\it Arabian Journal for Science
and Engineering}, {\bf34}(2), 1-11 (2009).


\bibitem{ChCZ} Yuqun Chen, Yongshan Chen and Chanyan Zhong, Composition-Diamond
Lemma for Modules, {\it Czech J. Math.}, to appear. arXiv:0804.0917


\bibitem{Ch}E.S. Chibrikov,  On free Lie conformal algebras, {\it Vestnik Novosibirsk
State University}, {\bf4}(1), 65-83 (2004).

\bibitem{Co} P.M. Cohn, On the embedding of rings in skew fields, {\it Proc. London
Math. Soc.},  {\bf11}(3), 511-530 (1961).


\bibitem{Ev} T. Evans, Embedding theorems for multiplicative systems and
projective geometries, {\it Proc. Amer. Math. Soc.}, {\bf3}, 614-620
(1952).


\bibitem{Fi1981} V.T. Filippov, On the chains of varieties generated by free Maltsev and alternative algebras,
{\it Dokl. Akad. Nauk SSSR}, {\bf 260}(5), 1082-1085 (1981).


\bibitem{Fi1984} V.T. Filippov, Varieties of Malcev and alternative algebras generated by algebras of finite
rank. Groups and other algebraic systems with finiteness conditions,
{\it Trudy Inst. Mat., Nauka Sibirsk. Otdel., Novosibirsk}, {\bf4},
139-156 (1984).

\bibitem{Go} A.P. Goryushkin, Imbedding of countable groups in
2-generated simple groups, {\it Math. Notes}, {\bf 16}, 725-727
(1974).


\bibitem{GuK08} L. Guo and W. Keigher,  On differential Rota-Baxter algebras,
 {\it Journal of Pure and Applied Algebra}, {\bf 212}, 522-540 (2008).


 \bibitem{Ha} P. Hall, On the embedding of a group in a join of given groups, {\it J.
Austral. Math. Soc.}, {\bf17}, 434-495 (1974).

 \bibitem{Hi} G. Higman, A finitely generated infinite simple group, {\it J. London Math. Soc.},
 {\bf26}, 61-64 (1951).

\bibitem{HNN49} G. Higman, B.H. Neumann and H. Neumann. Embedding theorems for
groups, {\it J. London Math. Soc.}, {\bf24}, 247-254 (1949).


\bibitem{Ho} John M. Howie, Fundamentals of Semigroup Theory, Clarendon
Press Oxford, 1995.


\bibitem{Iv65} I.S. Ivanov, Free T-sums of multioperator fields,
{\it Dokl. Akad. Nauk SSSR}, {\bf165}, 28-30 (1965); {\it Soviet
Math. Dokl.}, {\bf6}, 1398-1401 (1965).


\bibitem{Iv67} I.S. Ivanov, Free T-sums of multioperator fields, {\it Trudy Mosk.
Mat. Obshch.}, {\bf17}, 3-44 (1967).

 \bibitem{Ku}A.G. Kurosh, Free sums of multiple operator algebras,
{\it Siberian. Math. J.}, {\bf1}, 62-70 (1960).

 \bibitem{Ku2}A.G. Kurosh, Multioperators rings and algebras, {\it Russian Math.
Surveys. Uspekhi}, {\bf24}(1), 3-15 (1969).


\bibitem{Ly} R.C. Lyndon, On Burnside's problem, {\it Trans. Am. math. Soc.}, {\bf77}, 202-215 (1954).

\bibitem{LySc} R.C. Lyndon, P.E. Schupp, Combinatorial Group Theory,
Springer-Verlag, 1977.


 \bibitem{Ma} A.I. Malcev, On a representation of nonassociative rings [In Russian],
  {\it Uspekhi Mat. Nauk N.S.}, {\bf7},
 181-185 (1952).


\bibitem{ro}G.C. Rota, Baxter algebras and combinatorial identities I,
{\it Bull. Amer. Math. Soc.}, {\bf5}, 325-329 (1969).

\bibitem{She}I.P. Shestakov, A problem of Shirshov, {\it Algebra
Logic}, {\bf16}, 153-166 (1978).


\bibitem{Sh58}A.I. Shirshov, On free Lie rings, {\it Mat. Sb.}, {\bf45}, 13-21 (1958).

\bibitem{Sh56}A.I. Shirshov, On special J-rings, {\it Mat. Sbornik}, {\bf38}, 149-166 (1956).

\bibitem{Sh62}A.I. Shirshov, Some algorithmic problem for Lie algebras,
{\it Sibirsk. Mat. Z.}, {\bf3}, 292-296 (1962) (in Russian); English
translation in {\it SIGSAM Bull.}, {\bf33}(2), 3-6 (1999).


\bibitem{Sk} L.A. Skornyakov, Non-associative free T-sums of fields,
{\it Mat. Sb.}, {\bf44}, 297-312 (1958).


\bibitem{Shu} E.G. Shutov, Embeddings of semigroups in simple and complete
semigroups, {\it Mat. Sb.}, {\bf62}(4), 496-511 (1963).


\bibitem{Uf} V.A. Ufnarovski, combinatorial and Asymptotic Methods
in Algebra, {\it Encyclopaedia Mat. Sci.}, {\bf57}, 1-196 (1995).

\end{thebibliography}
\end{document}